\documentclass[12pt,leqno,draft]{article}

\usepackage{amssymb}
\usepackage[mathscr]{eucal}
\usepackage{amsmath,amssymb,latexsym,theorem,bbm}
\usepackage{color}
\usepackage{appendix}
\usepackage{setspace}
\usepackage[margin=1in]{geometry}
\usepackage[english]{babel}
\usepackage[round]{natbib}

\newcommand{\NN}{\mathbb{N}}

\newcommand{\RR}{\mathbb{R}}

\newcommand{\ZZ}{\mathbb{Z}}

\newcommand{\ta}{\widetilde{a}}

\newcommand{\tb}{\widetilde{b}}

\newcommand{\bd}{{\boldsymbol{d}}}

\newcommand{\tbd}{\widetilde{\bd}}

\newcommand{\bE}{{\boldsymbol{E}}}

\newcommand{\bI}{{\boldsymbol{I}}}

\newcommand{\bQ}{{\boldsymbol{Q}}}
\newcommand{\tbQ}{{\widetilde{\bQ}}}

\newcommand{\btheta}{{\boldsymbol{\theta}}}

\newcommand{\bone}{{\boldsymbol{1}}}

\newcommand{\cB}{{\mathcal B}}
\newcommand{\bcB}{\boldsymbol{\cB}}

\newcommand{\cD}{{\mathcal D}}

\newcommand{\cF}{{\mathcal F}}

\newcommand{\cL}{{\mathcal L}}
\newcommand{\cM}{{\mathcal M}}
\newcommand{\hM}{\widehat{M}}

\newcommand{\bcM}{\boldsymbol{\cM}}
\newcommand{\hbcM}{\widehat{\bcM}}

\newcommand{\cW}{{\mathcal W}}
\newcommand{\bcW}{\boldsymbol{\cW}}

\newcommand{\dd}{\mathrm{d}}
\newcommand{\ee}{\mathrm{e}}

\newcommand{\MLE}{{}}
\newcommand{\LSE}{{}}
\newcommand{\LSED}{{\mathrm{D}}}

\DeclareMathOperator*{\argmin}{arg\,min}

\newcommand{\EE}{\operatorname{\mathbb{E}}}
\newcommand{\PP}{\operatorname{\mathbb{P}}}
\newcommand{\OO}{\operatorname{O}}
\newcommand{\var}{\operatorname{Var}}
\newcommand{\cov}{\operatorname{Cov}}

\newcommand{\ha}{\widehat{a}}
\newcommand{\hb}{\widehat{b}}

\newcommand{\htau}{\widehat{\tau}}

\newcommand{\hbtheta}{\widehat{\btheta}}

\newcommand{\tbtheta}{\widetilde{\btheta}}

\newcommand{\tX}{\widetilde{X}}

\newcommand{\vare}{\varepsilon}

\renewcommand{\leq}{\leqslant}
\renewcommand{\geq}{\geqslant}

\newcommand{\stoch}{\stackrel{\PP}{\longrightarrow}}
\newcommand{\distr}{\stackrel{\cD}{\longrightarrow}}

\newcommand{\as}{\stackrel{{\mathrm{a.s.}}}{\longrightarrow}}

\newcommand{\proofend}{\hfill\mbox{$\Box$}}
\newcommand{\Norm}[1]{\left\lVert #1 \right\rVert}

\numberwithin{equation}{section}

\theoremstyle{change} \theorembodyfont{\em}
\newtheorem{Lem}{Lemma.}[section]
\newtheorem{Thm}[Lem]{Theorem.}
\newtheorem{Pro}[Lem]{Proposition.}
\newtheorem{Cor}[Lem]{Corollary.}

\theorembodyfont{\rm}
\newtheorem{Rem}[Lem]{Remark.}

\begin{document}

\begin{center}
 {\bfseries\Large Change detection in the Cox--Ingersoll--Ross model} \\[5mm]
 {\sc\large Gyula Pap} and  {\sc\large Tam\'as $\text{T. Szab\'o}^*$}
\end{center}

\vskip0.2cm

\noindent
 Bolyai Institute, University of Szeged,
 Aradi v\'ertan\'uk tere 1, H--6720 Szeged, Hungary.

\noindent e--mails: papgy@math.u-szeged.hu (G. Pap),
                    tszabo@math.u-szeged.hu (T. T. Szab\'o).

\noindent $*$ Corresponding author.

\renewcommand{\thefootnote}{}
\vspace*{0.2cm}
\footnote{This research was supported by the European Union and the State of Hungary, 
co-financed by the European Social Fund in the framework of T\'AMOP 4.2.4. 
A/2-11-1-2012-0001 ‘National Excellence Program’.
}

\vspace*{-10mm}

\begin{abstract}
We propose a change detection method for the famous
 Cox--Ingersoll--Ross model. This model is widely used in financial mathematics and therefore detecting a change in its parameters is of crucial importance. We develop one- and two-sided testing procedures for both drift parameters of the process. The test process is based on estimators that are motivated by the discrete time least-squares estimators, and its asymptotic distribution under the no-change hypothesis is that of a Brownian bridge. We prove the asymptotic weak consistence of the test, and derive the asymptotic properties of the change-point estimator under the alternative hypothesis of change at one point in time.
\end{abstract}

\section{Introduction}
\label{intro_pre}

We consider the well-known Cox--Ingersoll--Ross (CIR) model
 \begin{align}\label{CIR_SDE}
  \dd X_t = (a - b X_t) \, \dd t + \sigma \sqrt{X_t} \, \dd W_t , \qquad
  t \geq 0 ,
 \end{align}
 where \ $a > 0$, \ $b > 0$, \ $\sigma > 0$ \ and
 \ $(W_t)_{t \geq 0}$ \ is a standard Wiener process. We will be interested in detecting a change in the parameters $a$ and $b$, and for brevity we will use $\btheta:=(a,b)^{\top}.$ The volatility parameter $\sigma$ will not be estimated because we work with a continuous sample, from which (and indeed, from an arbitrarily small part of which) $\sigma$ can be calculated exactly, see \citet[Remark 2.6]{Barczy_13b}. Therefore change detection in $\sigma$ is not necessary -- we can calculate, without any uncertainty, whether $\sigma$ is constant across our sample. The constraints on the parameter values ensure the ergodic behavior of our process -- for details see Theorem \ref{Ergodicity} below. These constraints also ensure that any solution of \eqref{CIR_SDE} starting from a nonnegative value stays nonnegative indefinitely almost surely -- see Proposition \ref{Pro_CIR}.
 
 The process was proposed as an interest rate model by \citet{Cox_85} and is one of the standard "short rate" models in financial mathematics. The statistical properties of the model have therefore been extensively studied: \citet{Overbeck_98} provided estimators based on continuous-time observations, while the low-frequency discrete-time CLS estimators were proposed by \citet{Overbeck_97}. High-frequency estimators were proposed by \citet{Alaya_12, BenAlaya_13}, whose results we will require occasionally.
 
 There are a handful of change detection tests for the CIR process in the literature: \citet{Schmid_04} used control charts and a sequential method (i.e., an online procedure, which is in contrast to our offline one, where we assume the full sample to be known before starting investigations). They also supposed noisy observations, which will not be our interest. \citet{Guo_10} used the local parameter approach based on approximate maximum likelihood estimates. In essence, they wanted to find the largest interval for which the sample fits the model. Also, they used a discrete sample, whereas we will use a continuous one. The main result of our paper is that we were able to prove some asymptotic properties of the testing procedure under the alternative hypothesis as well as the null hypothesis. We believe this to be important because, if investigated only under the null hypothesis, a change-detection procedure is essentially a model-fitting test, and results under the alternative are necessary to verify its use for the more special task of change detection.
  
 The statistical problem we are concerned with is the following: we would like to test the null hypothesis
 \[
 \mathrm{H}_0: (X_t)_{t \in [0,T]} \text{ is the path of a CIR process}
 \]
 against the alternative hypothesis
  \begin{align*}
 \mathrm{H}_{\mathrm{A}}: \exists \tau \in [0,T]: & (X_{t})_{t\in[0,\tau]} \text{ is a CIR process with parameters } a=a', \ b=b', \text{ and } \\
 & (X_{t})_{t\in[\tau,T]} \text{ is a CIR process with parameters } a=a'', \ b=b''.
 \end{align*}
In general, we will be interested in asymptotic results as $T \to \infty.$ Under $\mathrm{H}_{\mathrm{A}}$ we will also require $\tau = \rho T$ with $\rho \in (0,1).$
 
 The layout of the paper is the following: in the remainder of the present Section 1 we will explain our notations. Section 2 will deal with the basic finite-sample and asymptotic properties of the CIR process and establishes the tools for our proofs. We will introduce our parameter estimators in Section 3 and derive their strong consistency. We will not investigate them in more detail than necessary since we will only use them to construct the test process, and we are more interested in their nice algebraic form than their statistical properties. We construct our test process and describe the test procedures in Section 4, where we also obtain the asymptotic distribution of the test process under $\mathrm{H}_{\mathrm{0}}$. Section 5 contains the first of our two asymptotic results -- namely, the weak consistence of the test. The second result, which concerns the properties of the change-point estimator under $\mathrm{H}_{\mathrm{A}}$, is stated and proved in Section 6. Section 7 explains how to modify the proofs in order to detect a change in $b$. Finally, the lemmata necessary for the proofs of the main theorems have been collected into Section 8.
 
\subsection{Notations}

In the following we describe our basic notations. Let \ $\NN$, \ $\ZZ_+$, \ $\RR$, \ $\RR_+$  \ and \ $\RR_{++}$ \ denote the sets
 of positive integers, non-negative integers, real numbers, non-negative real
 numbers and positive real numbers, respectively.
\ For \ $x , y \in \RR$, \ we will use
 \ $x \land y := \min(x, y)$ \ and \ $x \lor y := \max(x, y)$.
\ By \ $\|x\|$ \ and \ $\|A\|$ \ we denote the Euclidean norm of a vector
 \ $x \in \RR^d$ \ and the induced matrix norm of a matrix
 \ $A\in\RR^{d \times d}$, \ respectively.
We will use asymptotic notation for rates of convergence: $f(t) = \OO(g(t))$ means that $\limsup_{t \to \infty} \frac{f(t)}{g(t)} < \infty$. Similarly, for a stochastic process $X_t$, the notation $X_t = \OO_{\PP}(g(t))$ means that the collection of measures $\left(\cL\left(\frac{X_t}{g(t)}\right)\right)_{t \geq t_0}$ is tight for some $t_0 \in \RR_+$.
Unless otherwise noted, asymptotic statements are to be understood as $T \to \infty.$ 
Following the usual conventions, \ $\stoch$, \ $\distr$ \ and \ $\as$ \ will denote convergence
 in probability, in distribution and almost surely, respectively.

As for the probabilistic setup, \ $\bigl(\Omega, \cF, (\cF_t)_{t\in\RR_+}, \PP\bigr)$ \ will always be a filtered
 probability space satisfying the usual conditions, i.e.,
 \ $(\Omega, \cF, \PP)$ \ is complete, the filtration \ $(\cF_t)_{t\in\RR_+}$
 \ is right-continuous and \ $\cF_0$ \ contains all the $\PP$-null sets in
 \ $\cF$. We will repetadly work with continous martingales; as usual, their quadratic variation will be denoted by $\langle \cdot \rangle$.

\section{Preliminaries}
\label{Prel}

In our first proposition we recall some well-known properties of the solution of \eqref{CIR_SDE}.

\begin{Pro}\label{Pro_CIR}
For any random variable \ $\xi$ \ independent of \ $(W_t)_{t\in\RR_+}$
 \ and satisfying \ $\PP(\xi \in \RR_+) = 1$, \ there is a (pathwise) unique
 strong solution \ $(X_t)_{t\in\RR_+}$ \ of the SDE \eqref{CIR_SDE} with
 \ $X_0 = \xi$.
Further, we have \ $\PP(\text{$X_t \in \RR_+$ \ for all \ $t \in \RR_+$}) = 1$
 \ and the following equalities:
 \begin{align}\label{Solution}
  X_t &= \ee^{-bt}
        \left( X_0
               + a \int_0^t \ee^{bu} \, \dd u
               + \sigma \int_0^t \ee^{bu} \sqrt{X_u} \, \dd W_u \right) ,
  \qquad t \in \RR_+, \\ \label{eq:X2_solution}
  X_t^2 & = \ee^{-2bt} X_0^2 + \int_{0}^{t}\ee^{-2b(t-u)}(2a + \sigma^2)X_u \, \dd u + 2\sigma \int_{0}^{t} \ee^{-2b(t-u)} X_u^{3/2} \, \dd W_u, \qquad t \in \RR_+.
 \end{align}
 The conditional distribution of $X_t$ on $X_s$, where $s<t$, is noncentral chi-squared and we have
\begin{equation}\label{eq:finitesupmoments}
\sup_{t \in \RR_+} \EE(X_t^\eta) < \infty
\end{equation}
for all $\eta > 0$.
\end{Pro}

\noindent{\bf Proof.}
By a theorem due to Yamada and Watanabe \citep[see, e.g., ][Proposition 5.2.13]{Karatzas_91}, the strong uniqueness holds for
 \eqref{CIR_SDE}.
By \citet[Example V.8.2, page 221]{Ikeda_89}, there is a
 (pathwise) unique non-negative strong solution \ $(X_t)_{t\in\RR_+}$ \ of
 \eqref{CIR_SDE} with any initial value \ $\xi$ \
 independent of \ $(W_t)_{t\in\RR_+}$ \ and satisfying
 \ $\PP(\xi \in \RR_+) = 1$, \ and we have
 \ $\PP(\text{$X_t \in \RR_+$ \ for all \ $t \in \RR_+$}) = 1$.
\ Next, by application of the It\^o's formula for the process
 \ $(X_t)_{t\in\RR_+}$, \ we obtain
 \begin{align*}
  \dd(\ee^{bt} X_t)
  &= b \ee^{bt} X_t \, \dd t + \ee^{bt} \dd X_t
   = b \ee^{bt} X_t \, \dd t
     + \ee^{bt}
       \bigl( (a - b X_t) \, \dd t + \sigma \sqrt{X_t} \, \dd W_t \bigr) \\
  &= a \ee^{bt} \, \dd t + \sigma \ee^{bt} \sqrt{X_t} \, \dd W_t 
 \end{align*}
 for all \ $t \in \RR_+$, \ which implies \eqref{Solution}.  
 
 The noncentral chi-squared distribution is a well-known property of the process, and it can be found in the paper of \citet{Feller_51}. The property \eqref{eq:finitesupmoments} is a direct consequence of this fact and the calculations can be found, e.g., in \citet[Proposition 3]{BenAlaya_13}.
\proofend

The following result states the existence of a unique stationary distribution
 and the ergodicity of the CIR process. The proof can be put together from \citet{Feller_51}, \citet[Equation 20]{Cox_85}, and \citet{Jin_13}.

\begin{Thm}\label{Ergodicity}
Let \ $a, b, \sigma \in \RR_{++}$. 
\ Let \ $(X_t)_{t\in\RR_+}$ \ be a strong solution of \eqref{CIR_SDE} with
 \ $\PP(X_0 \in \RR_+) = 1$.
\ Then
 \renewcommand{\labelenumi}{{\rm(\roman{enumi})}}
 \begin{enumerate}
  \item
   $X_t \distr X_\infty$ \ as \ $t \to \infty$, \ and the distribution of
   \ $X_\infty$ \ is given by
   \begin{align}\label{Laplace}
    \EE(\ee^{-\lambda X_\infty})
    = \left(1 + \frac{\sigma^2}{2b} \lambda\right)^{-2a/\sigma^2} ,
    \qquad \lambda \in \RR_+ ,
   \end{align}
   i.e., \ $X_\infty$ \ has Gamma distribution with parameters \ $2a/\sigma^2$
   \ and \ $2b/\sigma^2$, \ hence 
   \[
     \EE(X_\infty^\alpha)
     = \frac{\Gamma\left(\frac{2a}{\sigma^2} + \alpha\right)}
            {\left(\frac{2b}{\sigma^2}\right)^\alpha
             \Gamma\left(\frac{2a}{\sigma^2}\right)} , \qquad
     \alpha \in \left(-\frac{2a}{\sigma^2}, \infty\right) .
   \]
 \item
  supposing that the random initial value \ $X_0$ \ has the same distribution
  as \ $X_\infty$, \ the process \ $(X_t)_{t\in\RR_+}$ \ is strictly stationary;
 \item
  for all Borel measurable functions \ $f : \RR \to \RR$ \ such that
  \ $\EE(|f(X_\infty)|) < \infty$, \ we have
  \begin{equation}\label{ergodic}
   \frac{1}{T} \int_0^T f(X_s) \, \dd s \as \EE(f(X_\infty)) \qquad
   \text{as \ $T \to \infty$.}
  \end{equation}
\end{enumerate}
\end{Thm}

\begin{Cor}
In the setting of Proposition \ref{Pro_CIR} we have
\begin{align*}
\EE(X_t) &= \ee^{-bt} \EE(X_0)
         + a \int_0^t \ee^{-b(t-u)} \, \dd u  \\
\EE(X_t^2) &= \ee^{-2bt} \EE(X_0^2) + \int_{0}^{t}(2a + \sigma^2)\left(\ee^{-b(2t-u)} \EE(X_0)
            + a \int_0^u \ee^{-b(2t-u-v)} \, \dd v\right) \, \dd u.
\end{align*}

Hence,
\begin{equation}\label{eq:X_mean_limits}
\lim_{t \to \infty} \EE(X_t) = \EE(X_\infty) = \frac{a}{b}, \qquad \lim_{t \to \infty} \EE(X_t^2) = \EE(X_\infty^2) = \frac{2a^2 + a^2 \sigma^2}{2b^2},
\end{equation}
moreover,
\begin{equation}\label{eq:X_diff_exp}
\int_0^{\infty} |\EE(X_t) - \EE(X_\infty)| \, \dd t < \infty, \qquad \int_0^{\infty} |\EE(X_t^2) - \EE(X_\infty^2)| \, \dd t < \infty.
\end{equation}
\end{Cor}

\noindent {\bf Proof.} The first equalities are straightforward by taking expectations on both sides in Proposition \ref{Pro_CIR} (we note that the stochastic integrals in question are indeed martingales due to \eqref{eq:finitesupmoments}). From there, \eqref{eq:X_mean_limits} is a question of elementary calculus: for the first equation we write
\begin{equation}\label{eq:EX_lim_calc}
\lim_{t \to \infty} \left(\ee^{-bt} \EE(X_0) + a \int_0^t \ee^{-b(t-u)} \, \dd u\right) = \lim_{t \to \infty} a \int_0^t \ee^{-bv} \, \dd v = a \int_0^{\infty} \ee^{-bv} \, \dd v = \frac{a}{b}.
\end{equation}
For the second equation we observe
\begin{equation}\label{eq:e_u_v_calc}
\int_{0}^{t}\int_{0}^{u} \ee^{-b(2t-u-v)} \, \dd v \dd u = \frac{1}{b}\left(\int_{0}^{t}(\ee^{-2b(t-u)} - \ee^{-b(2t-u)}) \, \dd u \right)= \frac{1}{b}\int_{0}^{t}\ee^{-2bu}\, \dd u + \frac{\ee^{-bt}}{b} \int_{0}^{t} \ee^{-bu} \, \dd u
\end{equation}
and hence
\begin{equation}\label{eq:EX2_lim_calc}
\begin{split}
\lim_{t \to \infty} & \left(\ee^{-2bt} \EE(X_0^2) + \int_{0}^{t}(2a + \sigma^2)\left(\ee^{-b(2t-u)} \EE(X_0)
+ a \int_0^u \ee^{-b(2t-u-v)} \, \dd v\right) \, \dd u \right) \\
& = (2a+ \sigma^2) \lim_{t \to \infty} \left( \EE(X_0) \ee^{-bt} \int_{0}^{t} \ee^{-bw} \, \dd w + a \int_0^t \int_0^u \ee^{-b(2t-u-v)} \, \dd v \, \dd u \right) \\
& = (2a+\sigma^2) \frac{1}{b}\int_0^{\infty} \ee^{-2bw} \, \dd w.
\end{split}
\end{equation}
 For the first part of \eqref{eq:X_diff_exp} we consider (keeping in mind \eqref{eq:EX_lim_calc})
\[
|\EE(X_t) - \EE(X_\infty)| = \left|\ee^{-bt} \EE(X_0) - a \int_{t}^{\infty} \ee^{-bu} \, \dd u \right| \leq \ee^{-bt} \EE(X_0) + ab^{-1} \ee^{-bt},
\]
which yields the result immediately. For the second part, we combine \eqref{eq:e_u_v_calc} and \eqref{eq:EX2_lim_calc} to obtain
\begin{align*}
|\EE(X_t^2) - \EE(X_\infty^2)| &= \left|\ee^{-2bt} \EE(X_0^2) + (2a+\sigma^2)\ee^{-bt}\int_{0}^{t}\left(\EE(X_0)\ee^{-bu} + \frac{1}{b} \ee^{-bu}\right) \, \dd u \right. \\
& \qquad \left. - \frac{1}{b} \int_{t}^{\infty}\ee^{-2bu}\, \dd u\right| \\
&\leq \ee^{-2bt}\EE(X_0^2) + (2a+\sigma^2)\ee^{-bt}\left(\EE(X_0) + \frac{1}{b}\right)\frac{1}{b} + \frac{1}{2b^2}\ee^{-2bt}.
\end{align*}
This yields the desired result immediately.
\proofend

Finally, we recall a strong law of large numbers and a central limit theorem for continuous local
 martingales.
\begin{Thm}{\bf \citep[Special case of][Lemma 17.4]{Liptser_01b}}
\label{DDS_stoch_int}
Let the process \ $(W_t)_{t\in\RR_+}$ \ be a standard Wiener process with respect to the
 filtration \ $(\cF_t)_{t\in\RR_+}.$
\ Let \ $(\xi_t)_{t\in\RR_+}$ \ be a measurable process adapted to $(\cF_t)_{t\in\RR_+}$ such that
\begin{align}\label{SEGED_STRONG_CONSISTENCY2}
   \PP\left( \int_0^t \xi_u^2 \, \dd u < \infty \right) = 1 ,
   \quad t \in \RR_+  \qquad \text{ and } \qquad \int_0^t \xi_u^2 \, \dd u \as \infty \qquad
   \text{as \ $t \to \infty.$}
\end{align} Then
 \begin{align}\label{SEGED_STOCH_INT_SLLN}
  \frac{\int_0^t \xi_u \, \dd W_u}
       {\int_0^t \xi_u^2 \, \dd u} \as 0 \qquad
  \text{as \ $t \to \infty$.}
 \end{align}
\end{Thm}

\begin{Thm}{\bf \citep[Special case of][Corollary VIII.3.24.]{Jacod_03}}\label{thm:MCLT}
Let $(X^n_t)_{t \in \RR_+}$ be a series of locally square-integrable continuous martingales such that
\[
\langle X^n \rangle_t \stoch t, \quad t \in \RR_+, \qquad \text{as } n \to \infty.
\]
Then $(X^n)_{t \in \RR_+} \distr (W_t)_{t \in \RR_+}$, where $(W_t)_{t \in \RR_+}$ is a standard Wiener process.
\end{Thm}

\section{Construction of our parameter estimators}
\label{section_EULSE}

In this section we will define some estimators for the drift parameters of the CIR process, based on continuous time observations. We will do this in the following way: first we introduce least squares estimators based on low-frequency discrete time observations, then we will introduce our estimators as a formal analogy; we will not try to construct our estimators as solutions to a least-squares problem.

An LSE of \ $(a, b)$ \ based on a discrete time observation
 \ $(X_i)_{i\in\{0,1,\ldots,n\}}$, \ can be obtained by solving the extremum
 problem
 \begin{align*}
  \bigl(\ha_n^{\LSED}, \hb_n^{\LSED}\bigr)
  := \argmin_{(a,b)\in\RR^2} \sum_{i=1}^n (X_i - X_{i-1} - (a - b X_{i-1}))^2 .
 \end{align*}
This is a simple exercise, which has the well-known solution
 \begin{align*}
  \begin{bmatrix}
   \ha_n^{\LSED} \\
   \hb_n^{\LSED}
  \end{bmatrix}
  = \begin{bmatrix}
     n & -\sum_{i=1}^n X_{i-1} \\
     -\sum_{i=1}^n X_{i-1} & \sum_{i=1}^n X_{i-1}^2
    \end{bmatrix}^{-1}
    \begin{bmatrix}
     \sum_{i=1}^n (X_i - X_{i-1})  \\
     - \sum_{i=1}^n (X_i - X_{i-1}) X_{i-1}
    \end{bmatrix} ,
 \end{align*}
 provided \ $n \sum_{i=1}^n X_{i-1}^2 - \left(\sum_{i=1}^n X_{i-1}\right)^2 > 0$.
 
By a formal analogy, we introduce the estimator of \ $(a, b)$ \ based on a continuous
 time observation \ $(X_t)_{t\in[0,T]}$ \ as
 \begin{align*}
  \hbtheta_T :=
  \begin{bmatrix}
   \ha_T^{\LSE} \\
   \hb_T^{\LSE}
  \end{bmatrix}
  &= \begin{bmatrix}
       T & -\int_0^T X_s \, \dd s \\
       -\int_0^T X_s \, \dd s & \int_0^T X_s^2 \, \dd s
      \end{bmatrix}^{-1}
      \begin{bmatrix}
       X_T - X_0 \\
      - \int_0^T X_s \, \dd X_s
      \end{bmatrix} ,
 \end{align*}
 provided
 \ $T \int_0^T X_s^2 \, \dd s - \left(\int_0^T X_s \, \dd s\right)^2 > 0$, which is true a.s. To see this, consider that, by a simple application of the Cauchy--Schwarz inequality,
 \[
T \int_0^T X_s^2 \, \dd s - \left(\int_0^T X_s \, \dd s\right)^2 \geq 0,
 \]
 and equality happens only if $X$ is constant almost everywhere on $[0,T]$. In particular, since $X$ is continuous on $[0,T]$ almost surely, this implies $X_0 = X_T$ almost surely. However, since the distribution of $X_T$ conditionally on $X_0 = x$ is absolutely continuous by Proposition \ref{Pro_CIR}, $\PP(X_T=X_0|X_0=x) = 0$ for all $x \in \RR_+$, which suffices for the statement.

To condense our notation, we will use
 \begin{equation}\label{def_Id}
\bQ_s:= \begin{bmatrix}
       s & -\int_0^s X_u \, \dd u \\
       -\int_0^T X_u \, \dd u & \int_0^s X_u^2 \, \dd u
      \end{bmatrix}
      \qquad \text{and} \qquad 
 \bd_s:=\begin{bmatrix}
      X_s - X_0       \\
     -\int_0^s X_u \, \dd X_u
       \end{bmatrix}
 \end{equation}

\begin{Rem}
The stochastic integral \ $\int_0^s X_u \, \dd X_u$ \ is observable, since,
 by It\^{o}'s formula, we have
 \ $\dd(X_t^2) = 2 X_t \,\dd X_t + \sigma^2 X_t \, \dd t$, \ $t \in \RR_+$,
 \ hence
 \[
   \int_0^s X_u \, \dd X_u
   = \frac{1}{2}
     \left( X_s^2 - X_0^2 - \sigma^2 \int_0^s X_u \, \dd u \right) .
 \]
\end{Rem}

Using the SDE \eqref{CIR_SDE} one can check that
 \begin{align}\label{eq:theta_difference}
   \begin{bmatrix}
   \ha_T^{\LSE} - a \\
   \hb_T^{\LSE} - b
  \end{bmatrix}
  &= \bQ_T^{-1}
     \begin{bmatrix}
      \sigma \int_0^T X_s^{1/2} \, \dd W_s \\
      - \sigma \int_0^T X_s^{3/2} \, \dd W_s
     \end{bmatrix} ,
 \end{align}
 provided
 \ $T \int_0^T X_s^2 \, \dd s - \left(\int_0^T X_s \, \dd s\right)^2 > 0$, which is, again, true a.s. In further calculations we will use
 \begin{equation}\label{def_tbd}
 \tbd_s := \sigma \begin{bmatrix}
       \int_0^s X_u^{1/2} \, \dd W_u \\
       - \int_0^s X_u^{3/2} \, \dd W_u
      \end{bmatrix}.
 \end{equation}

\begin{Thm}\label{Thm_LSE_cons}
\ Let \ $(X_t)_{t\in\RR_+}$ \ be a strong solution of \eqref{CIR_SDE} with
 \ $\PP(X_0 \in \RR_+) = 1$.
\ Then the LSE of \ $(a, b)$ \ is strongly consistent, i.e.,
 \ $\bigl(\ha_T^{\LSE}, \hb_T^{\LSE}\bigr) \as (a, b)$ \ as \ $T \to \infty$.
\end{Thm}

\noindent{\bf Proof.}
Recall \eqref{eq:theta_difference} and write
\begin{align}
   \begin{bmatrix}
   \ha_T^{\LSE} - a \\
   \hb_T^{\LSE} - b
  \end{bmatrix}
  &= \left(\frac{\bQ_T}{T}\right)^{-1}
     \frac{\sigma^2 \int_0^T X_s^3 \, \dd s}{T}
     \begin{bmatrix}
      \frac{\sigma^2 \int_0^T X_s \, \dd s}{\sigma^2 \int_0^T X_s^3 \, \dd s}\cdot \frac{\sigma \int_0^T X_s^{1/2} \, \dd W_s}{\sigma^2 \int_0^T X_s \, \dd s} \\
      \frac{- \sigma \int_0^T X_s^{3/2} \, \dd W_s}{\sigma^2 \int_0^T X_s^3 \, \dd s}
     \end{bmatrix}.
\end{align}
Now, the statement is evident from \eqref{SEGED_STOCH_INT_SLLN} and \eqref{ergodic}, noting that
\[
\frac{\sigma^2 \int_0^T X_s \, \dd s}{\sigma^2 \int_0^T X_s^3 \, \dd s} = \frac{T^{-1}\sigma^2 \int_0^T X_s \, \dd s}{T^{-1}\sigma^2 \int_0^T X_s^3 \, \dd s} \as \frac{\EE(X_\infty)}{\EE(X_\infty^3)}.
\]
\proofend

\section{Construction of the test process}
\label{section_test_process}

First we introduce the martingale
 \[
   M_s :=  X_s - X_0 - \int_0^s (a - b X_u) \, \dd u
        = \sigma \int_0^s \sqrt{X_u} \, \dd W_s , \qquad s \in \RR_+ ,
 \]
which satisfies
 \begin{equation}\label{eq:Mdiff}
   \dd M_s = \dd X_s - (a - b X_s) \, \dd s = \sigma \sqrt{X_u} \, \dd W_s .
 \end{equation}

Let us fix a time horizon \ $T \in \RR_{++}$.
\ The process will again be introduced as a formal analogy to the efficient score vector, as is done in \citet{Gombay_08}. The analogue of the efficient score vector process at time \ $tT$, \ $t \in [0, 1]$, \ will
 be 
 \begin{align*}
     \int_0^{tT}
      \begin{bmatrix}
       1 \\
       - X_s
      \end{bmatrix}
      \dd M_s .
 \end{align*}
The information contained in a continuous sample \ $(X_u)_{u\in[0,tT]}$ \ is the
 quadratic variation of the efficient score vector process, namely,
 \[
   \int_0^{tT}
    \begin{bmatrix}
     1 \\
     - X_s
    \end{bmatrix}
    \begin{bmatrix}
     1 \\
     - X_s
    \end{bmatrix}^\top
    \langle M \rangle_s \, \dd s
  = \sigma^2 \int_0^{tT}
    \begin{bmatrix}
     X_s & -X_s^2 \\
     -X_s^2  & X_s^3
    \end{bmatrix}
    \dd s
  =: \bI_{tT} ,
 \] 
 since \ $\langle M \rangle_s = \sigma^2 X_s$, \ $s \in \RR_+$.
\ For each \ $s \in \RR_+$, \ replacing the parameters by their estimates in
 \ $M_s$, \ we obtain an estimate \ $\hM_s^{(T)}$, \ i.e.,
 \[
   \hM_s^{(T)} := X_s - X_0- \int_0^s (\ha_T - \hb_T X_u) \, \dd u , \qquad
   s \in \RR_+ .
 \]
Our test process will be the estimated efficient score vector multiplied by
 the square root of the inverse of the information matrix, i.e.,
 \[
   \hbcM_t^{(T)}
   := \bI_T^{-1/2}
      \int_0^{tT}
      \begin{bmatrix}
       1 \\
       - X_s
      \end{bmatrix}
      \dd \hM_s^{(T)} , \qquad t \in[0, 1] .
 \]
This process can also be written in CUSUM form
 \[
   \hbcM_t^{(T)}
   = \bI_T^{-1/2} \bQ_{tT} \left( \hbtheta_{tT}^\MLE - \hbtheta_T^\MLE \right) ,
   \qquad t \in [0, 1] .
 \]
 Indeed,
 \begin{align*}
\int_0^{tT}
      \begin{bmatrix}
       1 \\
       - X_s
      \end{bmatrix}
      \dd \hM_s^{(T)} &= 
      \int_0^{tT}
      \begin{bmatrix}
       1 \\
       - X_s
      \end{bmatrix}
      \dd X_s -
\int_0^{tT}
      \begin{bmatrix}
       1 \\
       - X_s
      \end{bmatrix}
      \begin{bmatrix}
       1 \\
       - X_s
      \end{bmatrix}^{\top}
      \hbtheta_T
      \dd s
      \\ &=
      \bQ_{tT} \left(\bQ_{tT}^{-1}\int_0^{tT}
      \begin{bmatrix}
       1 \\
       - X_s
      \end{bmatrix}
      \dd X_s -\hbtheta_T \right).
 \end{align*}

\begin{Thm}\label{Thm_H0}
\ Let \ $(X_t)_{t\in\RR_+}$ \ be a strong solution of \eqref{CIR_SDE} with
 \ $\PP(X_0 \in \RR_+) = 1$.
\ Then
 \[
   \left(\hbcM_t^{(T)}\right)_{t\in[0,1]} \distr (\bcB_t)_{t\in[0,1]} \qquad
   \text{as \ $T \to \infty$,}
 \]
 where \ $(\bcB_t)_{t\in[0,1]}$ \ is a 2-dimensional standard Brownian bridge.
\end{Thm}

\noindent{\bf Proof.}
We have
 \begin{align*}
  \int_0^{tT}
   \begin{bmatrix}
    1 \\
    - X_s
   \end{bmatrix}
   \dd \hM_s^{(T)}
  &= \int_0^{tT}
      \begin{bmatrix}
       1 \\
       - X_s
      \end{bmatrix}
      \dd M_s
     - \int_0^{tT}
        \begin{bmatrix}
         1 \\
         - X_s
        \end{bmatrix}
        \left(\dd M_s - \dd \hM_s^{(T)}\right) ,
 \end{align*}
 and
 \begin{align*}
  \int_0^{tT}
   \begin{bmatrix}
    1 \\
    - X_s
   \end{bmatrix}
   \left(\dd M_s - \dd \hM_s^{(T)}\right)
  &= \int_0^{tT}
      \begin{bmatrix}
       1 \\
       - X_s
      \end{bmatrix}
      \left( \ha_T^\MLE - a - (\hb_T^\MLE - b) X_s \right)
      \dd s \\
  &= \int_0^{tT}
      \begin{bmatrix}
       1 \\
       - X_s
      \end{bmatrix}
      \begin{bmatrix}
       1 \\
       - X_s
      \end{bmatrix}^\top
      \begin{bmatrix}
       \ha_T^\MLE - a \\
       \hb_T^\MLE - b
      \end{bmatrix}
      \dd s
   = \bQ_{tT} \bQ_T^{-1} \tbd_T ,
 \end{align*}
 with the notations from \eqref{def_Id} and \eqref{def_tbd}. In the following, $\bE_2$ denotes the 2-dimensional identity matrix. From the preceding calculations it follows that, for every $t \in [0,1],$
 \begin{align*}
  \hbcM_t^{(T)}
  &= \bI_T^{-1/2} \left(\tbd_{tT} - \bQ_{tT} \bQ_T^{-1} \tbd_T\right) \\
  &= \bI_T^{-1/2} \left(\tbd_{tT} - t \tbd_T\right) + \bI_T^{-1/2}(t \bE_2-\bQ_{tT} \bQ_T^{-1})\tbd_T \\
  &= (T\bI)^{-1/2} \left(\tbd_{tT} - t \tbd_T\right) + ((T^{-1}\bI_T)^{-1/2} - \bI^{-1/2})T^{-1/2}\left(\tbd_{tT} - t \tbd_T\right) \\
  &\quad + \bI_T^{-1/2}(t \bE_2-\bQ_{tT} \bQ_T^{-1})\tbd_T,
 \end{align*}
  where
  \[
    \bI
    := \sigma^2 \begin{bmatrix}
    \EE(X_\infty) & -\EE(X_\infty^2) \\
    - \EE(X_\infty^2) & \EE(X_\infty^3)
    \end{bmatrix}.
  \]
 It is a simple consequence of the ergodic theorem that $T^{-1}\bI_T \as \bI$ as $T \to \infty$. 
 Consequently, Theorem \ref{Thm_H0} will follow from
 \begin{equation}\label{eq:informacio}
   \sup_{0 \leq t \leq 1}(t \bE_2-\bQ_{tT} \bQ_T^{-1}) \stoch 0 \qquad
   \text{as \ $T \to \infty$,}
   \end{equation}
 and
 \begin{equation}\label{eq:wiener}
   \left(T^{-1/2} \, \tbd_{tT}\right)_{t\in[0,1]}
   \distr (\bI^{1/2} \, \bcW_t)_{t\in[0,1]}
   \qquad \text{as \ $T \to \infty$,}
 \end{equation}
 where \ $(\bcW_t)_{t\in[0,1]}$ \ is a 2-dimensional standard Wiener process.

We begin by the proof of \eqref{eq:wiener}. The convergence is a simple consequence of Theorem \ref{thm:MCLT}. $\tbd_t$ is a locally square-integrable martingale, therefore we only need to check the pointwise convergence of the quadratic variation. Using (iii) from Theorem \ref{Ergodicity} it is easy to show that, for every $t \in [0,1]$,
\[
\frac{1}{T} \sigma^2 \int_0^{tT}
\begin{bmatrix}
X_s & -X_s^2 \\
 -X_s^2 & X_s^3
\end{bmatrix}
ds
\as
t
\begin{bmatrix}
\EE(X_\infty) & -\EE(X_\infty^2) \\
- \EE(X_\infty^2) & \EE(X_\infty^3)
\end{bmatrix}
=
t \bI,
\quad 
\text{as } T \to \infty.
\]
For \eqref{eq:informacio}, introduce
\[
    \bQ
    := \begin{bmatrix}
    1 & -\EE(X_\infty) \\
    - \EE(X_\infty) & \EE(X_\infty^2)
    \end{bmatrix}
  \]
  and note that due to Theorem \ref{Ergodicity} we have $T^{-1}\bQ_T \as \bQ.$
 Now, first observe that
\[
\lVert t \bE_2 - \bQ_{tT}\bQ_T^{-1}\rVert \leq t \left\lVert \frac{\bQ_T}{T} - \frac{\bQ_{tT}}{tT}\right\rVert \left\lVert\left(\frac{\bQ_T}{T}\right)^{-1}\right\rVert.
\]
For this transformation to be sensible, we needed to extend $\frac{\bQ_s}{s}$ continuously to $s=0$, but this can be done since all components of $\frac{I_s}{s}$ has a finite upper limit at 0 almost surely (i.e., the powers of $X_0$). Since the last factor converges to $\lVert\bQ^{-1}\rVert$ almost surely, for \eqref{eq:informacio} it is sufficient to show that
\begin{equation}\label{eq:suptconvergence}
\sup_{0 \leq t \leq 1} t \left\lVert \frac{\bQ_T}{T} - \frac{\bQ_{tT}}{tT}\right\rVert \stoch 0.
\end{equation}
In order to exploit the almost sure convergence of $\frac{\bQ_T}{T}$, we note that $\frac{\bQ_T}{T} \as \bQ$ implies $\sup_{s>T} \Norm{\frac{\bQ_s}{s}-\bQ} \as 0$ as $T \to \infty$ and thus $\sup_{s>T} \Norm{\frac{\bQ_s}{s}-\bQ} \stoch 0$ as $T \to \infty.$
Now let us introduce
\[
K:= \sup_{s \geq 0} \left \lVert \frac{\bQ_s}{s}\right \rVert.
\]
This limit is finite almost surely since $\frac{\bQ_s}{s}$ is continuous on $\RR_+$ and has a finite limit at infinity almost surely.
Now we observe, for an arbitrary $\epsilon > 0$,
\begin{align*}
\PP&\left(\sup_{0 \leq t \leq 1}t \left\lVert \frac{\bQ_T}{T} - \frac{\bQ_{tT}}{tT}\right\rVert > \epsilon\right) \\
& \leq \PP\left(\sup_{0 \leq t \leq \frac{\epsilon}{4 K} \wedge 1}t \left\lVert \frac{\bQ_T}{T} - \frac{\bQ_{tT}}{tT}\right\rVert > \epsilon\right) + \PP\left(\sup_{\frac{\epsilon}{4 K} \leq t \leq 1}t \left\lVert \frac{\bQ_T}{T} - \frac{\bQ_{tT}}{tT}\right\rVert > \epsilon\right) \\
& \leq \PP\left(\frac{\epsilon}{4 K}2K > \epsilon\right) + \PP\left(\sup_{\frac{\epsilon}{4 K} \leq t \leq 1} \left(t \left\lVert \frac{\bQ_T}{T} - \bQ \right\rVert + \left \lVert \frac{\bQ_{tT}}{tT} - \bQ \right\rVert \right) > \epsilon\right)\\
& \leq 0+\PP\left(\left\lVert \frac{\bQ_T}{T} - \bQ \right\rVert > \frac{\epsilon}{2}\right) + \PP\left(\sup_{\frac{\epsilon T}{4 K} \leq s} \left \lVert \frac{\bQ_{s}}{s} - \bQ \right\rVert > \frac{\epsilon}{2}\right).
\end{align*}
Dividing the last probability according to the value of $K$, we have
\begin{align*}
\PP&\left(\sup_{0 \leq t \leq 1}t \left\lVert \frac{\bQ_T}{T} - \frac{\bQ_{tT}}{tT}\right\rVert > \epsilon\right) \\
& \leq \PP\left(\left\lVert \frac{\bQ_T}{T} - \bQ \right\rVert > \frac{\epsilon}{2}\right) + \PP\left(\left\{\sup_{\frac{\epsilon T}{4 K} \leq s} \left \lVert \frac{\bQ_{s}}{s} - \bQ \right\rVert > \frac{\epsilon}{2}\right\} \bigcap \left\{K \leq \sqrt{T}\right\}\right) \\
& \quad + \PP (K > \sqrt{T}) \\
& \leq \PP\left(\left\lVert \frac{\bQ_T}{T} - \bQ \right\rVert > \frac{\epsilon}{2}\right) + \PP\left(\sup_{\frac{\epsilon \sqrt{T}}{4} \leq s} \left \lVert \frac{\bQ_{s}}{s} - \bQ \right\rVert > \frac{\epsilon}{2}\right) + \PP \left(K > \sqrt{T}\right).
\end{align*}
All three terms in the last expression tend to zero as $T \to \infty$, therefore \eqref{eq:informacio} is proved.
\proofend

\subsection{Testing procedures}\label{subsec:testing}

Based on Theorem \ref{Thm_H0}, we can develop the following tests with a significance level of $\alpha$:

\noindent {\bf Test 1 (one-sided):} if it is clear that, in case of a change, $a' < a''$, reject $\mathrm{H}_0$ if the minimum of $(\hbcM^{(1)}_t)_{t \in [0,T]}$ is greater than $C_1(\alpha)$, where $C_1(\alpha)$ can be obtained from the distribution of the minimum of a standard Brownian bridge. The same test can be applied to the maximum (for $a' > a''$) and to $(\hbcM^{(2)}_t)_{t \in [0,T]}$ (for a change in $b$).

\noindent {\bf Test 2 (two-sided):} reject $\mathrm{H}_0$ if the maximum of $|\hbcM^{(1)}_t|_{t \in [0,T]}$ is greater than $C_2(\alpha)$, where $C_2(\alpha)$ can be obtained from the distribution of the maximum of the absolute value of standard Brownian bridge. The same test can be applied to $|\hbcM^{(2)}_t|_{t \in [0,T]}$ (for a change in $b$).

\noindent Naturally, the test for $a$ and $b$ can be applied simultaneously, in which case the significance levels for the individual tests have to be modified accordingly, in order to produce an overall significance level of $\alpha$.

\section{Asymptotic consistence of the test}

Before stating our results under the alternative hypothesis, we need to examine the ergodicity results that we can use more closely. Let us take two parameter vectors: $\btheta'$ and $\btheta''$ (in the formulation of the theorem, $\btheta' = (a',b)^{\top}$ and $\btheta'' = (a'', b)^{\top}$, but for the time being, we can work more generally). Furthermore, we take two random variables, $X_\infty'$ and $X_\infty''$, such that they are distributed according to the stationary distributions corresponding to $\btheta'$ and $\btheta''$, respectively. Let us take a process $(X_t)_{t \in \RR_+}$ such that it evolves according to \eqref{CIR_SDE} with parameters $\btheta'$ until $t=\rho T$ and with parameters $\btheta''$ thereafter. We would like to apply the ergodic theorem (i.e., Theorem \ref{Ergodicity}) separately to the process before and after the change-point (i.e., $\rho T$). However, we cannot do this directly for the second part because the initial distribution may depend on $T$. However, we do have
\begin{equation}\label{aergodic}
\frac{1}{T - \rho T}\int_{\rho T}^T g(X_t) \, \dd t \stoch \EE(g(\tX'')),
\end{equation}
where \ $g:\RR_+ \to \RR$ \ with \ $\EE(|g(\tX)|) < \infty$. 
\ Indeed, for an arbitrary \ $\vare>0$ 
\begin{align*}
&\PP\left( \left| \frac{1}{T - \rho T}\int_{\rho T}^T g(X_t) \, \dd t
- \EE(g(\tX'')) \right| > \vare \right) \\
&= \int_{\RR_+}
\PP\left( \left| \frac{1}{T - \rho T}\int_{\rho T}^T g(X_t) \, \dd t
- \EE(g(\tX''))\right| > \vare
\, \Bigg| \, X_{\rho T}=x \right)
\, \dd P_X^{\rho T}(x) \\
& \leq \left \lVert P_X^{\rho T} - P^* \right \rVert  + \int_{\RR_+}
\PP\left( \left| \frac{1}{T - \rho T}\int_{\rho T}^T g(X_t) \, \dd t
- \EE(g(\tX''))\right| > \vare
\, \Bigg| \, X_{\rho T}=x \right)
\, \dd P^*(x),
\end{align*}
where $P^*$ is the distribution of $\tX'$, $P_X^{\rho T}$ is the distribution of $X_{\rho T}$ and $\lVert \cdot \rVert$ is the total variation norm. The first term converges to zero because the CIR process is positive Harris recurrent \citep[][Theorem 2.5]{Jin_13}. This implies ergodicity by \citet[Theorem 6.1]{Meyn_93}, since in this case the 1-skeleton (i.e., the process $(X_i)_{i \in \ZZ_+}$) is clearly irreducible because the support of the distribution of $X_1$ conditionally on $X_0$ is $\RR_+$. In the second term the measure is finite, while the integrand is bounded by 1 and converges to zero pointwise, therefore \eqref{aergodic} is proved by the Lebesgue Dominated Convergence Theorem. The same line of reasoning can be used to apply Theorem \ref{DDS_stoch_int} (with weak convergence) and Theorem \ref{thm:MCLT} after the point of change.
Let us now introduce
\[
\bd_{[a,b]}:= \begin{bmatrix}
       \int_a^b
        1 \dd X_s \\[2mm]
       - \int_a^b X_s \dd X_s
      \end{bmatrix} , \qquad
   \bQ_{[a,b]}
   := \begin{bmatrix}
       \int_a^b 1 \, \dd s
        & -\int_a^b X_s \, \dd s \\[2mm]
       -\int_a^b  X_s \, \dd s
        & \int_a^b X_s^2 \, \dd s
      \end{bmatrix}.
\]
With these notations,
\[
\hbtheta_T^{\MLE} = \left(\bQ_{[0,\tau]}+\bQ_{[\tau,T]}\right)^{-1}(\bd_{[0,\tau]}+\bd_{[\tau,T]}).
\]
With the help of the ergodic theorem, we can see that this quantity has a finite weak limit:
\[
\tbtheta := \begin{bmatrix} \ta\\ \tb \end{bmatrix} :=
 (\rho \bQ' + (1-\rho) \bQ'')^{-1} \left(\rho \bQ' \btheta'
+ (1-\rho) \bQ''
\btheta'' \right),
\]
where
\[
\bQ' := \begin{bmatrix}
1 & - \EE(X_\infty') \\
- \EE(X_\infty') & \EE((X_\infty')^2)
\end{bmatrix}, \qquad
\bQ'' := \begin{bmatrix}
1 & - \EE(X_\infty'') \\
- \EE(X_\infty'') & \EE((X_\infty'')^2)
\end{bmatrix}.
\]

\begin{Thm}\label{thm:consistence}
If $a$ changes from $a'>0$ to $a''>0$ at time $\tau = \rho T$, where $\rho \in (0,1)$, then for any $\gamma \in \left(0, \frac{1}{4}\right)$ we have
\[
\sup_{0 \leq t \leq T} \hM_t^{(T)} = T \psi + \OO_{\PP}(T^{1-\gamma}),
\]
with $\psi = (a' - a'') \bone_1^{\top} ( (\rho \bQ')^{-1} + ((1-\rho)\bQ'')^{-1})^{-1} \bone_1.$ Here $\bone_1 = (1,0)^{\top}$, the first unit vector.  
\end{Thm}

\begin{Rem} Note how the sign of the principal term depends on the direction of change: it is negative in case of an upwards change and positive in case of a downwards change. This gives us the possibility to design one-sided tests.
\end{Rem}

\noindent {\bf Proof}. First we show how the estimates behave in this case.  
Clearly,
\[
\btheta' - \tbtheta = (1-\rho) (\rho \bQ' + (1-\rho) \bQ'')^{-1} \bQ'' (\btheta' - \btheta'').
\]
We have
\[
\begin{bmatrix}
1 \\ -\EE(X_\infty')
\end{bmatrix}
=
\bQ' \begin{bmatrix}
1 \\ 0
\end{bmatrix},
\]
and hence
\begin{equation}\label{eq:theta'_psi}
(\btheta' - \tbtheta)^{\top} \begin{bmatrix}
1 \\ -\EE(X_\infty')
\end{bmatrix} = \frac{\psi}{\rho}.
\end{equation}
In the same way we can conclude that
\begin{equation}\label{eq:theta''_psi}
(\btheta'' - \tbtheta)^{\top} \begin{bmatrix}
1 \\ -\EE(X_\infty'')
\end{bmatrix} = -\frac{\psi}{1-\rho}.
\end{equation}
Now we apply the following decomposition (which is useful for $t< \tau$; for $t \geq \tau$ it has to be modified in a straightforward manner):
\begin{equation}\label{eq:decomp}
\begin{split}
\int_0^{t}1 \dd \hM_u^{(T)} &= \int_{0}^{t}1\dd M_u + \int_{0}^{t}\left[(\tb-b')\EE(X_u) + (a'-\ta)\right]\dd u \\
& \quad + \int_{0}^{t} (\hb_T-b') (X_u - \EE(X_u)) \, \dd u + \int_{0}^{t} \left[(\hb_T-\tb) \EE(X_u) + (\ta - \ha_T)\right] \dd u \\
& = \int_{0}^{t}1\dd M_u + \int_{0}^{t} (\btheta' - \tbtheta)^{\top} \begin{bmatrix} 1 \\ - \EE(X_u) \end{bmatrix} \, \dd u \\
& \quad + \int_{0}^{t} (\btheta' - \hbtheta_T) \begin{bmatrix} 0 \\ \EE(X_u) - X_u \end{bmatrix} \, \dd u +
\int_{0}^{t} (\tbtheta - \hbtheta_T) \begin{bmatrix} 1 \\ - \EE(X_u) \end{bmatrix} \, \dd u.
\end{split}
\end{equation}
This leads to
\begin{align*}
& \left| \sup_{0 \leq t \leq T} \int_0^{t}1 \dd \hM_u^{(T)} - T \psi \right| \\
& \leq \sup_{0 \leq t \leq T} \left| \int_{0}^{\tau \wedge t}1\dd M_u' + \int_{\tau \wedge t}^{t}1\dd M_u'' \right|
  + \sup_{0 \leq t \leq T} \left| \int_{0}^{t} (\tbtheta - \hbtheta_T) \begin{bmatrix} 1 \\ - \EE(X_u) \end{bmatrix} \, \dd u \right| \\
& \quad + \left|\sup_{0 \leq t \leq  T} \left(\int_{0}^{\tau \wedge t} (\btheta' - \tbtheta)^{\top} \begin{bmatrix} 1 \\ - \EE(X_u) \end{bmatrix} \, \dd u  + \int_{\tau \wedge t}^{t} (\btheta'' - \tbtheta)^{\top} \begin{bmatrix} 1 \\ - \EE(X_u) \end{bmatrix} \, \dd u - T \psi\right)\right| \\
& \quad + \sup_{0 \leq t \leq T} \left| \int_{0}^{\tau \wedge t} (\btheta' - \hbtheta)^{\top} \begin{bmatrix} 0 \\ \EE(X_u)- X_u \end{bmatrix} \, \dd u + \int_{\tau \wedge t}^{t} (\btheta'' - \hbtheta)^{\top} \begin{bmatrix} 1 \\ \EE(X_u) - X_u \end{bmatrix} \, \dd u \right|.
\end{align*}
The first term is $\OO_{\PP}(T^{\gamma - 1})$ according to Lemma \ref{lem:M_sup}, the fourth term by Lemma \ref{lem:X_EX_diff_sup} and the second term by Lemma \ref{lem:theta_conv_alt}. For the third term, we write
\begin{align*}
& \left|\sup_{0 \leq t \leq  T} \left(\int_{0}^{\tau \wedge t} (\btheta' - \tbtheta)^{\top} \begin{bmatrix} 1 \\ - \EE(X_u) \end{bmatrix} \, \dd u + \int_{\tau \wedge t}^{t} (\btheta'' - \tbtheta)^{\top} \begin{bmatrix} 1 \\ - \EE(X_u) \end{bmatrix} \, \dd u - T \psi \right) \right| \\
&\leq \sup_{0 \leq t \leq  T} \left|\int_{0}^{\tau \wedge t} (\btheta' - \tbtheta)^{\top} \begin{bmatrix} 0 \\ \EE(X_\infty) - \EE(X_u) \end{bmatrix} \, \dd u \right| \\
& \quad + \sup_{0 \leq t \leq  T} \left| \int_{\tau \wedge t}^{t} (\btheta'' - \tbtheta)^{\top} \begin{bmatrix} 1 \\ \EE(X_\infty) - \EE(X_u) \end{bmatrix} \, \dd u\right| + \left| \sup_{0 \leq t \leq  T} \left(\frac{\tau \wedge t}{\rho} - \frac{(t-\tau)^+}{1-\rho} - T \right) \psi \right|.
\end{align*}
The first two terms in this decomposition are bounded by \eqref{eq:X_diff_exp} and the last one is obviously zero, with the supremum attained at $t=\tau$. This completes the proof. \proofend

\section{Estimation of the change point}

The natural estimate of the change point if $a' > a''$, i.e., when a downward change in $a$ is being tested, is
\[
\htau_T := \inf \{t \in \RR_+: \hM_t^{(T)} = \sup_{0 \leq t \leq T} \hM_t^{(T)} \}.
\]
Clearly, this is a well-defined, finite quantity, since $\hM_t^{(T)}$ has continuous trajectories almost surely. Regarding this estimate, we state the following result:
\begin{Thm}\label{thm:changepoint}
Under the asumptions of Theorem \ref{thm:consistence}, if $a' > a''$, then we have
\[
\htau_T - \rho T = \OO_{\PP}(1).
\]
\end{Thm}

\noindent {\bf Proof.} 
We remind the reader that, according to the assumptions, $\tau = \rho T$. We need to show that
\[
\lim_{K \to \infty} \sup_{T \in \RR}\PP(|\htau_T - \rho T| \geq K) = 0 \quad \text{a.s.},
\]
or, equivalently,
\[
\lim_{K \to \infty} \limsup_{T \in \RR}\PP(|\htau_T - \rho T| \geq K) = 0 \quad \text{a.s.}
\]

For this, it is sufficient to show that
\begin{equation}\label{eq:max_before}
\lim_{K \to \infty} \limsup_{T \to \infty} \PP \left( \sup_{\rho T - K < t < \rho T + K} \hM_t^{(T)} \leq \sup_{0 \leq t \leq \rho T - K} \hM_t^{(T)} \right) = 0
\end{equation}
and that
\begin{equation}\label{eq:max_after}
\lim_{K \to \infty} \limsup_{T \to \infty} \PP \left( \sup_{\rho T - K < t < \rho T + K} \hM_t^{(T)} \leq \sup_{\rho T + K \leq t \leq T} \hM_t^{(T)} \right) = 0.
\end{equation}
First we prove \eqref{eq:max_before}. We observe
\begin{align*}
&\PP\left(\sup_{\rho T - K < t < \rho T + K} \hM_t^{(T)} \leq \sup_{0 \leq t \leq \rho T - K} \hM_t^{(T)} \right)
 \leq \PP\left(\hM_{\rho T}^{(T)} \leq \sup_{0 \leq t \leq \rho T - K} \hM_t^{(T)}\right) \\
& = \PP \left(\inf_{0 \leq t \leq \rho T - K} (\hM_{\rho T}^{(T)} - \hM_t^{(T)}) \leq 0\right)
  = \PP\left(\inf_{K \leq t \leq \rho T} t^{-1} \int_{\rho T - t}^{\rho T}  1 \, \dd \hM_s^{(T)} \leq 0\right)
\end{align*}

We apply the decomposition \eqref{eq:decomp} to show that
\begin{equation}\label{eq:tau_decomp}
\begin{split}
\PP & \left(\inf_{K \leq t \leq \rho T} t^{-1} \int_{\rho T - t}^{\rho T}  1 \, \dd \hM_s^{(T)} \leq 0\right) \\
& \leq \PP \left( \inf_{K \leq t \leq \rho T} t^{-1} \int_{\rho T - t}^{\rho T} (\btheta' - \tbtheta)^{\top} \begin{bmatrix}
1 \\ - \EE(X_s)
\end{bmatrix}
\, \dd s \leq \frac{\psi}{2}\right) \\
& \quad + \PP\left(\sup_{K \leq t \leq \rho T}\left|t^{-1}(M_{\rho T} - M_{\rho T - t})\right|\geq \frac{\psi}{6}\right) \\
& \quad + \PP \left(\sup_{K \leq t \leq \rho T} \left|t^{-1} \int_{\rho T-t}^{\rho T} (\btheta' - \hbtheta_T)^{\top} 
\begin{bmatrix}
0 \\ \EE(X_s) - X_s
\end{bmatrix}
\, \dd s \right| \geq \frac{\psi}{6}\right) \\
& \quad + \PP \left(\sup_{K \leq t \leq \rho T} \left|t^{-1} \int_{\rho T - t}^{\rho T} (\tbtheta - \hbtheta_T)^{\top}
\begin{bmatrix}
1 \\ - \EE(X_s)
\end{bmatrix}
\, \dd s \right| \geq \frac{\psi}{6}\right).
\end{split}
\end{equation}
In the first term we take the probability of a deterministic event, therefore it is either 0 or 1; we show that for sufficiently large $K,N$ it is 0.

Actually, this is the same statement in continuous time as Lemma 7.7 in \citet{Pap_13}, and the proof is also the same. First we note that, as has been shown before,
\[
f(t):= (\btheta' - \tbtheta)^{\top} \begin{bmatrix}
1 \\ - \EE(X_t)
\end{bmatrix} \to \psi, \quad t \to \infty.
\]
For an arbitrary $\vare>0$, let us introduce $\nu(\vare):=\sup_{t: f(t) < \psi - \vare} < \infty$. Furthermore, let 
\[
\kappa(\vare) := \inf_{0 \leq t \leq \nu(\vare)} f(t) > -\infty. 
\]
Then we have, for a sufficiently large $T$,
\begin{align*}
&\inf_{0 \leq t \leq \rho T} t^{-1} \int_{\rho T -t}^{\rho T} f(s) \, \dd s \\
& \geq \min\left(\inf_{0 \leq t \leq \rho T-\nu(\vare)} t^{-1} \int_{\rho T -t}^{\rho T} f(s) \, \dd s, \frac{1}{\rho T - \nu(\vare)} \left(\kappa(\vare)\nu(\vare) + (\rho T-\nu(\vare)) \inf_{\nu(\vare) \leq s \leq \rho T} f(s)\right)\right) \\
& \geq \min\left(\psi - \vare, \frac{\kappa(\vare)\nu(\vare)}{\rho T - \nu(\vare)} + \psi - \vare\right).
\end{align*}
As $\frac{\kappa(\vare)\nu(\vare)}{\rho T - \nu(\vare)} \to 0$ as $T \to \infty$, we conclude that the second term in \eqref{eq:tau_decomp} is 0 for sufficiently large $T$, irrespective of $K$.

The fourth term in \eqref{eq:tau_decomp} converges to zero as $T \to \infty$ for any $K$, as
\[
\sup_{0 \leq t \leq \rho T} \left\lVert t^{-1} \int_{\rho T -t}^{\rho T} \begin{bmatrix}
1 \\ -\EE(X_s)
\end{bmatrix}
\, \dd s \right\rVert \leq \sup_{0 \leq t \leq \rho T} \left\lVert 
\begin{bmatrix}
1 \\ -\EE(X_t)
\end{bmatrix} \right\rVert,
\]
and the right hand side is bounded as $T \to \infty$. Meanwhile, $\tbtheta - \hbtheta_T \to 0$ a.s., which suffices for the fourth term in \eqref{eq:tau_decomp}.
For the third term we use Lemma \ref{lem:lim_limsup_X} and for the second one we can use Lemma \ref{lem:M_sup}. \proofend

\section{Detecting a change in $b$}

In Theorems \ref{thm:consistence} and \ref{thm:changepoint} we postulated a change in $a$. However, this was only done to keep the resulting calculations tractable. Straightforward modifications allow us to prove the same results for a change in $b$ -- following the same thoughts as in \ref{subsec:testing}. In this case, we would have
\begin{Thm}\label{thm:b_consistence}
If $b$ changes from $b'>0$ to $b''>0$ at time $\tau = \rho T$, where $\rho \in (0,1)$, then for any $\gamma \in \left(0, \frac{1}{4}\right)$ we have
\[
\sup_{0 \leq t \leq T} \int_0^t (- X_s) \dd \hM_s^{(T)} = T \phi + \OO_{\PP}(T^{1-\gamma}),
\]
with $\phi = (b' - b'') \bone_2^{\top} (((1-\rho)\bQ'')^{-1} + (\rho \bQ')^{-1})^{-1} \bone_2.$ Here $\bone_2 = (0,1)^{\top}$, the second unit vector.  
\end{Thm}

In place of \eqref{eq:decomp} we have, then,
\begin{equation}
\begin{split}
\int_0^{t} X_u \dd \hM_u^{(T)} &= \int_{0}^{t} X_u \dd M_u + \int_{0}^{t}\left[(\tb-b')\EE(X_u^2) + (a'-\ta)\EE(X_u)\right]\dd u \\
& \quad + \int_{0}^{t} \left[ (\hb_T-b') (X_u^2 - \EE(X_u^2)) + (a'-\ta)(X_u -\EE(X_u)) \right] \, \dd u \\
& \quad + \int_{0}^{t} \left[(\hb_T-\tb) \EE(X_u^2) + (\ta - \ha_T)\EE(X_u)\right] \dd u \\
& = \int_{0}^{t}X_u\dd M_u + \int_{0}^{t} (\btheta' - \tbtheta)^{\top} \begin{bmatrix} -\EE(X_u) \\ \EE(X_u^2) \end{bmatrix} \, \dd u \\
& \quad + \int_{0}^{t} (\btheta' - \hbtheta_T) \begin{bmatrix} \EE(X_u)-X_u \\ X_u^2 - \EE(X_u^2) \end{bmatrix} \, \dd u +
\int_{0}^{t} (\tbtheta - \hbtheta_T) \begin{bmatrix} -\EE(X_u) \\ \EE(X_u^2) \end{bmatrix} \, \dd u.
\end{split}
\end{equation}
The second term is approximately
\[
t (\btheta' - \tbtheta)^{\top} \begin{bmatrix} -\EE(X_\infty') \\ \EE((X_\infty'')^2) \end{bmatrix} = t (\btheta' - \tbtheta)^{\top} \bQ' \begin{bmatrix} 0 \\ 1 \end{bmatrix} = t \frac{\phi}{\rho}
\]
as in \eqref{eq:theta'_psi}. From here the proofs proceed as for a change in $a$, with the only added difficulty that we will require \eqref{eq:X2int_var} and as well as \eqref{eq:Xint_var}, but the proof of that result is merely a matter of algebra.

\section{Details of the proofs}\label{sec:lemmata}

In this section we detail the necessary lemmata for the proofs of our main theorems. Some of them, especially Lemma \ref{lem:Xint_var}, are rather technical and depend essentially on tedious but straightforward calculations. Others, while using more sophisticated tools, are also tailored to the specific needs of the proofs and their proofs are not particularly insightful themselves, hence they were relegated to this section. The one exception to this is Lemma \ref{lem:KoLe}, which is an analogue of Lemma \ref{kole} and may deserve independent interest.

\begin{Lem}\label{lem:Xint_var}
For the model described by \eqref{CIR_SDE} we have
\begin{equation}\label{eq:Xint_var}
\var \left(\int_{0}^{t}X_s \dd s \right) = \OO (t), \quad t \to \infty,
\end{equation}
and
\begin{equation}\label{eq:X2int_var}
\var \left(\int_{0}^{t}X_s^2 \dd s \right) = \OO (t), \quad t \to \infty.
\end{equation}
\end{Lem}

\noindent {\bf Proof.}
For \eqref{eq:Xint_var} we note
\[
\var \left(\int_{0}^{t}X_s \dd s \right) = \EE \left(\int_{0}^{t}(X_u - \EE X_u) \dd u \int_{0}^{t}(X_v - \EE X_v) \dd v\right) = \iint\limits_{[0,t]^2} \cov (X_u,X_v) \dd u \dd v.
\]
By using \eqref{Solution}, we can write
\begin{equation}\label{eq:CovXuXv}
\begin{split}
\cov(X_u, X_v) &= \EE [(X_u - \EE X_u)(X_v - \EE X_v)] = \\
&= \EE \left[\left(e^{-bu}(X_0-\EE X_0) + \sigma \int_{0}^{u}e^{-b(u-w)}\sqrt{X_w} \dd W_w\right) \right. \\
& \quad
\left. \times \cdot \left(e^{-bv}(X_0-\EE X_0) + \sigma \int_{0}^{v}e^{-b(v-z)}\sqrt{X_z} \dd W_z\right)\right] \\
&= e^{-b(u+v)}\var(X_0) + \sigma^2 \int_{0}^{u \wedge v} e^{-b(u+v-2w)}\EE(X_w)\dd w \\
&\leq e^{-b(u+v)}\var(X_0) + (E(X_0) + \frac{a}{b})\sigma^2 \int_{0}^{u \wedge v} e^{-b(u+v-2w)} \dd w,
\end{split}
\end{equation}
since
\[
\EE(X_w) = e^{-bw}\EE(X_0) + a \int_{0}^{w}e^{-bs} \dd s
\]
by \eqref{Solution}. Furthermore,
\begin{equation}\label{eq:e_uv_int_est}
\begin{split}
&\iint\limits_{[0,t]^2} \left(\int_{0}^{u \wedge v} e^{-b(u+v-2w)} \dd w\right) \dd u \dd v = \iint\limits_{[0,t]^2} \left[\frac{e^{-b(u+v-2w)}}{2b}\right]_{w=0}^{w=u \wedge v} \dd u \dd v \\
& \quad = \iint\limits_{[0,t]^2} \left[\frac{1}{2b}\left(e^{-b|u-v|}-e^{-b(u+v)}\right)\right] \dd u \dd v \leq \frac{1}{b}\iint\limits_{[0,t]^2} e^{-b|u-v|}\dd u \dd v = \OO(t).
\end{split}
\end{equation}
We combine this with the last line of \eqref{eq:CovXuXv} and note that
\begin{equation}\label{eq:e_uv_int_asymp}
\iint\limits_{[0,t]^2} e^{-b(u+v)} \dd u \dd v = \OO(t),
\end{equation}
which completes the proof of \eqref{lem:Xint_var}.

For \eqref{eq:X2int_var} we use the same approach. By \eqref{eq:X2_solution},
\begin{align*}
\cov(X_u^2, X_v^2) &= \EE[(X_u^2 - \EE(X_u^2))(X_v^2 - \EE(X_v^2))] = \\
& = \EE \left[ \left( \ee^{-2bu} (X_0^2 - \EE(X_0)^2)  + \int_{0}^{u}(2 a+\sigma^2) \ee^{-b(2u-w)} (X_0 - \EE(X_0)) \, \dd w \right. \right. \\
& \qquad \qquad + (2a+\sigma^2)\sigma \int_{0}^{u}\ee^{-2b(u - w)} \int_{0}^{w} \ee^{-b(w-z)} \sqrt{X_z} \, \dd W_z \\
& \qquad \qquad \left. + \sigma \int_{0}^{u}\ee^{-2b(u-w)} X_w^{3/2} \, \dd W_w \right) \\
& \qquad \times \left( \ee^{-2bv} (X_0^2 - \EE(X_0)^2)  + \int_{0}^{v}(2 a+\sigma^2) \ee^{-b(2v-w)} (X_0 - \EE(X_0)) \, \dd w \right. \\
& \qquad \qquad + (2a+\sigma^2)\sigma \int_{0}^{v}\ee^{-2b(v - w)} \int_{0}^{w} \ee^{-b(w-z)} \sqrt{X_z} \, \dd W_z \\
& \qquad \qquad \left. \left. + \sigma \int_{0}^{u}\ee^{-2b(v-w)} X_w^{3/2} \, \dd W_w \right) \right] \\
& = \ee^{-2b(u+v)} \var(X_0^2) + (2a+\sigma^2)^2\ee^{-b(u+v)}\int_0^u \ee^{-bw} \, \dd w \int_0^v \ee^{-bw} \, \dd w \var(X_0) \\
& \quad + 2 \ee^{-3b(u+v)}(2a+\sigma^2)\int_0^u \ee^{-bw} \, \dd w \int_0^v \ee^{-bw} \, \dd w \cov(X_0 ,X_0^2) \\
& \quad + (2a+\sigma^2)^2 \sigma^2 \int_{0}^{u} \int_{0}^{v} \EE \left( \int_{0}^{w} \ee^{-b(2u-w-z)} \sqrt{X_z} \dd W_z \right. \\
& \quad \phantom{+ (2a+\sigma^2)^2 \sigma^2 \int_{0}^{u} \int_{0}^{v} \EE \Bigg(}\times \left. \int_{0}^{r} \ee^{-b(2v-r-q)} \sqrt{X_q}\dd W_q\right)\dd r \dd w \\
& \quad + (2a+\sigma^2) \sigma^2 \int_{0}^{u} \EE \left( \int_{0}^{w} \ee^{-b(2u-w-z)} \sqrt{X_z} \dd W_z \int_{0}^{v} \ee^{-2b(v-q)} X_q^{3/2}\dd W_q\right) \dd w \\
& \quad + (2a+\sigma^2) \sigma^2 \int_{0}^{v} \EE \left( \int_{0}^{w} \ee^{-b(2v-w-z)} \sqrt{X_z} \dd W_z \int_{0}^{u} \ee^{-2b(u-q)} X_q^{3/2}\dd W_q\right) \dd w \\
& \quad + \int_{0}^{u \wedge v} e^{-2b(u+v-2w)}\EE(X_w^3)\dd w.
\end{align*}
Proceeding with the calculations, we have
\begin{align*}
\cov(X_u^2, X_v^2) & \leq \ee^{-2b(u+v)} \var(X_0^2) + (2a+\sigma^2)^2 \frac{1}{b^2}\ee^{-b(u+v)} \var(X_0) \\
& \quad + 2 \ee^{-3b(u+v)}(2a+\sigma^2)\frac{1}{b^2} \cov(X_0 ,X_0^2) \\
& \quad + (2a+\sigma^2)^2 \sigma^2 \int_{0}^{u} \int_{0}^{v} \int_{0}^{w \wedge r} \ee^{-b(2u+2v-w-r-2z)} \EE(X_z) \dd z \dd r \dd w \\
& \quad + (2a+\sigma^2) \sigma^2 \int_{0}^{u} \int_{0}^{w \wedge v} \ee^{-b(2u+2v-w-3z)} \EE(X_z^2) \dd z \dd w \\
& \quad + (2a+\sigma^2) \sigma^2 \int_{0}^{v} \int_{0}^{w \wedge u} \ee^{-b(2u+2v-w-3z)} \EE(X_z^2) \dd z \dd w \\
& \quad + \int_{0}^{u \wedge v} e^{-2b(u+v-2w)}\EE(X_w^3)\dd w.
\end{align*}
Referring to \eqref{eq:e_uv_int_est} and \eqref{eq:e_uv_int_asymp} we see that we need only concern ourselves about the fourth, fifth and sixth terms. For the fifth term, we have, for $u < v$,
\begin{align*}
\int_{0}^{u} \int_{0}^{w \wedge v} \ee^{-b(2u+2v-w-3z)} \dd z \dd w &= \int_{0}^{u} \int_{0}^{w} \ee^{-b(2u+2v-w-3z)} \dd z \dd w \\
&= \int_0^u \frac{1}{3b}(\ee^{-b(2u+2v-4w)} - \ee^{-b(2u+2v-w)}) \dd w \\
&= \frac{1}{12b^2}(\ee^{-2b(v-u)} - \ee^{-2b(v+u)}) - \frac{1}{3b^2} (\ee^{-b(u+2v)} - \ee^{-2b(u+v)}),
\end{align*}
and for $u > v$,
\begin{align*}
\int_{0}^{u} \int_{0}^{w \wedge v} \ee^{-b(2u+2v-w-3z)} \dd z \dd w &= \int_{0}^{v} \int_{0}^{w} \ee^{-b(2u+2v-w-3z)} \dd z \dd w \\
& \quad + \int_v^u \int_0^v \ee^{-b(2u+2v-w-3z)} \dd z \dd w \\
& = \int_0^v \frac{1}{3b}(\ee^{-b(2u+2v-4w)} - \ee^{-b(2u+2v-w)}) \dd w \\
& \quad + \int_v^u \frac{1}{3b}(\ee^{-b(2u-v-w)} - \ee^{-b(2u+2v-w)}) \dd w \\
& = \frac{1}{12b^2} \left( \ee^{-2b(u-v)} - \ee^{-2b(u+v)}\right) - \frac{1}{3b^2}\left(\ee^{-b(2u+v)} - \ee^{-2b(u+v)}\right) \\
& \quad + \frac{1}{3b^2}\left( \ee^{-b(u-v)} - \ee^{-2b(u-v)} - \left( \ee^{-b(u+2v)} - \ee^{-b(2u+v)} \right) \right).
\end{align*}
The same results, with $u$ and $v$ exchanged, hold for the sixth term. All the exponential expressions in question can be estimated from above by $\ee^{-b|u-v|},$ whence we can invoke \eqref{eq:e_uv_int_est} again to conclude that the fifth and sixth terms, integrated over $[0,t]^2$, are $\OO(t).$

All that remains is the fourth term: for $u < v$,
\begin{align*}
\int_{0}^v & \int_{0}^{v} \int_{0}^{w \wedge r} \ee^{-b(2u+2v-w-r-2z)} \dd z \dd r \dd w \\
&= \int_{0}^{u} \int_{0}^{w} \int_{0}^{r} \ee^{-b(2u+2v-w-r-2z)} \dd z \dd r \dd w  + \int_{0}^{u} \int_{w}^{r} \int_{0}^{r} \ee^{-b(2u+2v-w-r-2z)} \dd z \dd r \dd w \\
& = \int_{0}^u \int_{0}^{w} \frac{1}{2b} \left( \ee^{-b(2u+2v-w-3r)} - \ee^{-b(2u+2v-w-r)} \right) \dd r \dd w \\
& \quad + \int_{0}^{u} \int_{w}^{v} \frac{1}{2b} \left( \ee^{-b(2u+2v-3w-r)} - \ee^{-b(2u+2v-w-r)}\right) \dd r \dd w \\
& = \int_{0}^{u} \left[\frac{1}{6b^2}\left(\ee^{-b(2u+2v-4w)} - \ee^{-b(2u+2v-w)}\right) - \frac{1}{2b^2} \left(\ee^{-b(2u+2v-2w)} - \ee^{-b(2u+2v-w)}\right) \right. \\
& \phantom{= \int_{0}^{u} \Bigg[}\left. + \frac{1}{2b^2} \left(\ee^{-b(2u+v-3w)} - \ee^{-b(2u+2v-4w)}\right) - \frac{1}{2b^2} \left(\ee^{-b(2u + v - w)} - \ee^{-b(2u+2v-2w)}\right) \right] \dd w \\
& = \frac{1}{24b^3}\left(\ee^{-2b(v-u)} - \ee^{-2b(u+v)}\right) - \frac{1}{6b^3} \left(\ee^{-b(u+2v)} - \ee^{-2b(u+v)}\right) \\
& \quad - \frac{1}{4b^3}\left(\ee^{-2bv} - \ee^{-2b(u+v)}\right) + \frac{1}{2b^3} \left(\ee^{-b(u+2v)} - \ee^{-b(2u+2v)}\right) \\
& \quad + \frac{1}{6b^3}\left(\ee^{-b(v-u)} - \ee^{-b(2u+v)}\right) - \frac{1}{8b^3}\left(\ee^{-2b(v-u)} - \ee^{-2b(u+v)}\right) \\
& \quad - \frac{1}{2b^3}\left(\ee^{-b(u+v)} - \ee^{-b(2u+v)}\right) + \frac{1}{4b^3}\left(\ee^{-2bv} - \ee^{-2b(u+v)}\right),
\end{align*}
and, $u$ and $v$ have to be interchanged for $u > v$ (in this case, we exchange the two outer integrals, and from there, the modifications are trivial). Again, we see that all the exponential terms are dominated by $\ee^{-b|u-v|}$, which, by invoking \eqref{eq:e_uv_int_est}, completes the proof of the lemma, noting that $\sup_{t \geq 0} \EE(X_t^3) < \infty.$
\proofend

The following lemma is an analogue of Lemma \ref{kole}, which is a H\'ajek--R\'enyi type inequality. With Lemma \ref{kole} one can estimate the tail probabilities of the maximum of a random sequence, based solely on the joint moments of the elements and, critically, without the assumption of independence. In our applications, not the supremum of a sequence but the maximum of a function is considered, so we had to modify the statement accordingly.

It turns out that the proof can be constructed along the lines of Theorem 4.1 in \citet{Kokoszka_00}. In that paper, a slightly stronger result than Lemma \ref{kole} was formulated and proven; however, it was impractical to use, hence the more useful corollary formulated as Theorem 3.1 in \citet{Kokoszka_98}, which is obtainable from Theorem 4.1 in \citet{Kokoszka_00} by a simple application of the Cauchy--Schwarz theorem.

\begin{Lem}\label{lem:KoLe}
Let $Y_t$ be a process with a.s. continuous trajectory, $\alpha, \beta \in \RR_+$ with $\alpha < \beta$ and $c$ a deterministic function. Then, for any $\vare > 0$,
\begin{align*}
\vare^2\PP&\left\{ \sup_{s \in [\alpha,\beta]} \left(c(s) \int_{0}^{s}Y_u \dd u \right)^2 > \vare^2 \right\} \leq c(\alpha)^{2} \int_0^{\alpha} \EE(Y_u^2) \, \dd u \\ & + \int_{\alpha}^{\beta}\left(\int_{0}^{s}\int_{0}^{s}\EE(Y_uY_v)\dd u \dd v\right) \dd |c(s)^2| + 2 \int_{\alpha}^{\beta} c(s)^2\left[\EE(Y_s^2)\int_{0}^{s}\int_0^s\EE(Y_uY_v)\dd u \dd v\right]^{1/2} \dd s
\end{align*}
\end{Lem}

\noindent {\bf Proof.}
For any nonnegative process $Z_t$ with a.s. continuous trajectories and a.s. locally bounded variation, let $\tau_{\vare}$ be the first hitting time of $[\vare, \infty)$ in $[\alpha, \infty)$, $A$ be the event $\{\tau_{\vare} < \beta\}$ and $D_s$ be the event $\{\sup_{\alpha \leq u \leq s} Z_u \leq \vare\}.$ Note that $D_\beta = A^C$. Then it is easy to check that
\begin{equation}\label{eq:epsilon_Z}
\vare \bone_A \leq Z_\alpha + \int_{\alpha}^{ \beta}\bone_{D_s} \dd Z_s.
\end{equation}
Indeed, if $A$ occurs, the LHS is $\vare$, and the RHS is $\vare$, if $Z_\alpha < \vare$ and $Z_\alpha$ if $Z_\alpha \geq \vare$. If $A^C$ occurs, the LHS is zero, while the RHS is $Z_\beta \geq 0$.

Let us apply this result with $Z_t = c(t)^2\left\lvert\int_{0}^{t} Y_s \, \dd s \right\rvert^2$. We take expectations on both sides:
\begin{align*}
\vare^2& \PP\left(\sup_{\alpha \leq s \leq \beta} \left \lvert c(s) \int_{0}^{s}Y_u \dd u \right \rvert > \vare \right) \\
& \leq \EE\left[ c(\alpha)^{2} \int_{0}^{\alpha} Y_u^2 \, \dd u \right] + \EE \left[\int_{\alpha}^{\beta}\bone_{D_s} \dd \left(\left(c(s)\int_{0}^{s}Y_u \dd u\right)^2\right)\right] \\
& = c(\alpha)^{2} \int_0^{\alpha} \EE(Y_u^2) \, \dd u + \EE \left[2 \int_{\alpha}^{\beta} \bone_{D_s} c(s) \int_{0}^{s} Y_u \dd u \left(\left(\int_{0}^{s}Y_u \dd u \right) \dd c(s) + c(s) Y_s \dd s\right)\right] \\
& = c(\alpha)^{2} \int_0^{\alpha} \EE(Y_u^2) \, \dd u \\
& \quad + \EE \left[2 \int_{\alpha}^{\beta} \bone_{D_s} \left( \int_{0}^{s} \int_{0}^{s} Y_u Y_v \dd u \dd v \right) \dd (c^2(s)) + 2 \int_{\alpha}^\beta \bone_{D_s} c^2(s) Y_s \int_{0}^{s} Y_u \dd u \dd s\right] \\
& \leq c(\alpha)^{2} \int_0^{\alpha} \EE(Y_u^2) \, \dd u \\
& \quad + \EE \left[2 \int_{\alpha}^{\beta} \bone_{D_s} \left( \int_{0}^{s} \int_{0}^{s} Y_u Y_v \dd u \dd v \right) \dd \lvert c^2(s)\rvert + 2 \int_{\alpha}^\beta \bone_{D_s} c^2(s) Y_s \int_{0}^{s} Y_u \dd u \dd s\right] 
\end{align*}
In the last step we replaced the induced norm of $c^2(s)$ by its total variation norm. Indeed, the inequality holds because $\int_{0}^{s} \int_{0}^{s} Y_u Y_v \dd u \dd v = \left(\int_{0}^{s}Y_u \dd u\right)^2$ for every $\omega$ in the probability space where $Y$ is defined, therefore the integrand is nonnegative. Now, we employ several well-known inequalities and the replacement of the indicator function by 1 to obtain our statement. \proofend

\begin{Lem}\label{lem:X_EX_diff_sup}
If the parameters $a$ and $b$ remain constant, we have, for any $\gamma < \frac{1}{4}$,
\[
\sup_{0 \leq t \leq T} t^{\gamma - 1} \int_{0}^{t}|X_u - \EE(X_u)| \, \dd u = \OO_{\PP}(1).
\]
\end{Lem}

\noindent {\bf Proof.} We will use Lemma \ref{lem:KoLe} for the process $Y_t:= X_t - \EE(X_t)$ and $c(s) = s^{\gamma - 1}$ and $\alpha=0, \ \beta = T$. Then we can use Lemma \ref{lem:Xint_var} to conclude that
\[
\int_{0}^{s} \int_{0}^{s} \EE(Y_uY_v) \, \dd v \, \dd u = \int_{0}^{s} \int_{0}^{s} \cov(X_u, X_v) \, \dd v \, \dd u \leq \kappa s, \quad s \in \RR_+,
\]
for some constant $\kappa > 0$. Hence, in this case,
\begin{align*}
\int_{0}^{T}&\left(\int_{0}^{s}\int_{0}^{s}\EE(Y_uY_v)\dd u \dd v\right) \dd |c(s)^2| + 2 \int_{0}^{T} c(s)^2\left[\EE(Y_s^2)\int_{0}^{s}\int_0^s\EE(Y_uY_v)\dd u \dd v\right]^{1/2} \dd s \\
& \leq \int_{0}^{T} \kappa (2 - 2 \gamma) s^{2 \gamma - 2} \, \dd s + 2 \int_0^T s^{2 \gamma -2} (K \kappa s)^{1/2} \, \dd s \\
& = \kappa (2-2\gamma)\int_0^T s^{2\gamma - 2} \, \dd s + 2 (K \kappa)^{1/2} \int_0^T s^{2 \gamma - 3/2} \, \dd s < \infty.
\end{align*}
This implies the desired statement immediately. \proofend

\begin{Lem}\label{lem:M_sup}
If the parameters $a$ and $b$ remain constant, we have, for any $\gamma < \frac{1}{2}$,
\[
\sup_{0 \leq t \leq T} T^{\gamma - 1} |M_t| = \OO_{\PP}(1).
\]
\end{Lem}

\noindent {\bf Proof.}
First we note that $(M_t)_{t \in \RR_+}$ has an a.s. continuous trajectory on $\RR_+$, therefore also on $[0,1]$. Thus we conclude that $\sup_{0 \leq t \leq 1} |M_t| = \OO_{\PP(1)}$. Next, we use the law of the iterated logarithm for continuous martingales. This can be put together from the Dambis--Dubins--Schwarz theorem \citep[Theorem 3.4.6]{Karatzas_91} and the law of the iterated logarithm for the Wiener process \citep[Theorem 2.9.23]{Karatzas_91}.
\[
\limsup_{t \to \infty} \frac{|M_t|}{\sigma^{2\lambda}\left(\int_{0}^{t}X_u \, \dd u\right)^{\lambda}} \leq \limsup_{t \to \infty} \frac{|M_t|}{\sigma \sqrt{\int_{0}^{t}X_u \, \dd u} \sqrt{\log \log (\sigma^2 \int_{0}^{t}X_u \, \dd u)}} = 1 \quad \text{a.s.}, \quad \forall \lambda > \frac{1}{2},
\]
which means that the supremum on $[1, \infty]$ is finite a.s. (since the process in question has a.s. continuous trajectories). Now we note that
\[
\frac{\sigma^{2\lambda}\left(\int_{0}^{t}X_u \, \dd u\right)^{\lambda}}{t^\lambda} \to \sigma^{2\lambda} \EE(X_\infty)^{\lambda} \quad \text{a.s.}.
\]
Now the statement of the lemma is obtained straightforwardly since
\[
\sup_{0 \leq t \leq T} T^{\gamma - 1} |M_t| = \max(\sup_{0 \leq t \leq 1} T^{\gamma - 1} |M_t|, \sup_{1 \leq t \leq T} T^{\gamma - 1} |M_t|),
\]
and both terms have been shown to be $\OO_{\PP}(1)$. \proofend

\begin{Lem}\label{lem:theta_conv_alt}
Under the conditions of Theorem \ref{thm:consistence} we have
\[
\hbtheta - \tbtheta = \OO_{\PP} (T^{-1/2}).
\]
\end{Lem}

\noindent {\bf Proof.}
\begin{align*}
T^{1/2}(\hbtheta - \tbtheta) & = (T^{-1} \bQ_T)^{-1} T^{-1/2} \left[\bd_{0,\tau}-\bQ_T \tbQ^{-1} \left(\rho \bQ' \begin{bmatrix}
a \\ b
\end{bmatrix} \right) \right. \\ & \phantom{(T^{-1} \bQ_T)^{-1} T^{-1/2}} \qquad \left. +
\bd_{\tau,T} - \bQ_T \tbQ^{-1} \left((1-\rho) \bQ'' \begin{bmatrix}
a'' \\ b''
\end{bmatrix}\right) \right]
\end{align*}

The first factor converges almost surely, so we analyze
\begin{align*}
T^{-1/2}&\left[\bd_{0,\tau}-\bQ_T \tbQ^{-1} \left(\rho \bQ' \begin{bmatrix}
a \\ b
\end{bmatrix} \right)\right] \\
&= T^{-1/2} \tbd_\tau + T^{-1/2}\left( \bQ_{[0,\tau]} - \bQ_T \tbQ^{-1} \rho \bQ'
\right)
\begin{bmatrix}
a' \\ b'
\end{bmatrix}  
\end{align*}
The first term is $\OO_{\PP}(1)$ by \eqref{eq:wiener}. We need to show that the second term is also $\OO_{\PP}(1)$. For this, we can neglect the vector of the parameters, which are constant, so we investigate
\begin{align*}
T^{-1/2}\left( \bQ_\tau - \rho \bQ_T \tbQ^{-1} \bQ'\right) &= T^{-1/2} \left( \bQ_\tau - \EE(\bQ_{\tau})\right)  + T^{-1/2} \left(\EE(\bQ_{\tau}) - \tau \bQ'\right) \\ &\quad  - T^{-1/2} \left(\rho (\bQ_T - \EE(\bQ_T)) \tbQ^{-1} \bQ'\right) \\
& \quad - T^{-1/2} \left( \rho (\EE(\bQ_T) - T \tbQ) \tbQ^{-1} \bQ' \right)
\end{align*}
The first and third factors have a finite variance at the limit, by Lemma \ref{lem:Xint_var}. Therefore, by an application of Chebyshev's inequality, we have that they are $\OO_{\PP}(1)$. The second and fourth terms are deterministic and $\OO(1)$ by \eqref{eq:X_diff_exp}. \proofend

\begin{Lem}\label{lem:lim_limsup_X}
Under the conditions of Theorem \ref{thm:consistence} we have
\[
\lim_{K \to \infty} \limsup_{T \to \infty} \PP \left( \sup_{K \leq t \leq \rho T} \left| t^{-1} \int_{\rho T - t}^{\rho T} (X_s - \EE(X_s)) \, \dd s \right| > \frac{\psi}{6} \right) = 0.
\]
\end{Lem}
\noindent {\bf Proof.}
We use Lemma \ref{lem:KoLe}. We choose $c(s) = s^{-1}$ and $Y_s = X_{\rho T - s} - \EE(X_{\rho T - s})$ with $\alpha = K$ and $\beta = \rho T$. The estimate on the probability in question is then
\begin{equation}\label{eq:prob_est_1}
\begin{split}
K^{-2} &\int_{\rho T - K}^{\rho T} \var(X_u) \, \dd u + \int_{K}^{\rho T}\left(\int_{\rho T - s}^{\rho T}\int_{\rho T -s }^{\rho T}\cov(X_u,X_v)\dd u \dd v\right) \dd \left|s^{-2}\right| \\
& \quad + 2 \int_{K}^{\rho T} s^{-2}\left[\var(X_s)\int_{\rho T -s}^{\rho T}\int_{\rho T-s}^{\rho T}\cov(X_u,X_v)\dd u \dd v\right]^{1/2} \dd s.
\end{split}
\end{equation}
Now we make use of \eqref{eq:CovXuXv} and \eqref{eq:e_uv_int_est} to show that
\begin{align*}
\int_{\rho T -s}^{\rho T}\int_{\rho T-s}^{\rho T}\cov(X_u,X_v)\dd u \dd v & \leq \var(X_0) \int_{\rho T -s}^{\rho T}\int_{\rho T-s}^{\rho T} \ee^{-b(u+v)} \dd u \dd v \\
& \quad + (\EE(X_0) + ab^{-1})\sigma^2 b^{-1} \int_{\rho T -s}^{\rho T}\int_{\rho T-s}^{\rho T} \ee^{-b|u-v|} \dd u \dd v
 \leq \mu s,
\end{align*}
for some positive constant $\mu$.
We introduce $\lambda:= \sup_{t \in \RR} \var(X_t) < \infty,$  to continue the estimation started in \eqref{eq:prob_est_1}:
\[
K^{-2} K \lambda + 2 \int_K^\rho T s^{-3} \mu s \, \dd s + 2 \int_K^{\rho T} s^{-2} (\lambda \mu)^{1/2} s^{1/2} \, \dd s.
\]
Clearly, as $T \to \infty$ (and hence $\rho T \to \infty$), and then $K \to \infty$, this expression tends to zero, which completes our proof. \proofend

\begin{Lem}
Under the conditions of Theorem \ref{thm:consistence} we have, for any $\vare > 0$,
\[
\lim_{K \to \infty} \limsup_{T \to \infty} \PP \left( \sup_{K \leq t \leq \rho T} \left| t^{-1} (M_\rho T - M_{\rho T-t}) \right| > \vare \right) = 0.
\]
\end{Lem}

\noindent {\bf Proof.} Let us take a backward partition of $[0, \rho T]$ such that $0 = t_n < t_{n-1} < t_{n-2} < \ldots < t_1 < t_0 = \rho T.$ For $t \in [t_{i+1}, t_i]$, we have
\[
\left|\frac{M_\rho T - M_t}{\rho T -t}\right| \leq \left|\frac{M_\rho T - M_{t_{i+1}}}{\rho T - t_i}\right| + \left| \frac{M_t - M_{t_{i+1}}}{\rho T - t_i} \right|.
\]
Therefore, we have the following estimation:
\begin{equation}\label{eq:disc_cont_decomp}
\begin{split}
& \PP \left( \sup_{K \leq t \leq \rho T} \left| t^{-1} (M_\rho T - M_{\rho T-t}) \right| > \vare \right) = \PP \left( \sup_{0 \leq t \leq \rho T-K} \left| (\rho T-t)^{-1} (M_\rho T - M_t) \right| > \vare \right) \\
& \leq \PP\left( \max_{i^* \leq i \leq n} \left| (\rho T - t_{i})^{-1} (M_\rho T - M_{t_{i+1}}) \right| > \frac{\vare}{2} \right) \\
& \quad + \sum_{i=i*}^{n} \PP\left(\sup_{t_{i+1} < t < t_{i}} \left| (\rho T - t_{i})^{-1} (M_t - M_{t_{i+1}}) \right| > \frac{\vare}{2} \right),
\end{split}
\end{equation}
where $i^* = \min\{i: t_i < \rho T -K\}.$
Let us use this estimate with $t_i := \rho T - 2^{i-1}$ for $0 < i < n$, so that $n = \lfloor \log_2 \rho T \rfloor$ and $i^* = \lfloor \log_2 K \rfloor + 1$.

For the first term we can use the following lemma:
 \begin{Lem} {\bf \citep[Theorem 3.1]{Kokoszka_98}} \label{kole}
Let \ $(Y_n)_{n\in\NN}$ \ be a sequence of random variables with finite second moments, and let \ $(c_n)_{n \in \NN}$ \ be a sequence of nonnegative constants. Then, for any \ $a > 0$,
 \begin{align*}
  a^2 \PP \left( \max_{1 \leq k \leq n} c_k \left| \sum_{j=1}^k Y_j \right| 	> a \right)
  &\leq \sum_{k=1}^{n-1} |c^2_{k+1} - c^2_k| \sum_{i,j=1}^k \EE(Y_iY_j) \\
  &\quad +2 \sum_{k=1}^{n-1} c_{k+1}^2 \left(\EE\left( Y^2_{k+1} \right) \sum_{i,j=1}^k \EE\left(Y_iY_j\right)\right)^{1/2} \\
  &\quad +2 \sum_{k=0}^{n-1} c_{k+1}^2 \EE(Y_{k+1}^2).
 \end{align*}
\end{Lem}
Let us set $Y_1:=M_{t_{i^*}}-M_{\rho T}$, and $Y_k=M_{t_{i^*+k-1}} -M_{t_{i^*+k-2}}$ for $1 < k \leq n-i^*+1$ and $c_k = (\rho T - t_{i^* + k - 1})^{-1}$. Let us note that due to the structure of the $t_i$, we have $c_k = 2^{-(i^* + k - 2)}$ for $k \leq n- i^* $ and $2^{-(n-1)} < c_{n-i^*+1} < c_{n-i^*}$. Consequently, we can use
\[
|c_{k+1}^2 - c_k^2| \leq |c_{k+1}^2 - 4 c_{k+1}^2| = 3 c_{k+1}^2.
\]
Also, notice that
\begin{align*}
\sum_{i,j=1}^k \EE(Y_iY_j) &= \EE\left(\sum_{i=1}^k Y_i\right)^2 = \EE(M_{t_{i^*+k-1}} - M_{\rho T})^2 \\
&= \sigma^2 \int_{t_{i^*+k-1}}^{\rho T}\EE(X_u) \, \dd u \leq \sigma^2 \mu (\rho T - t_{i^*+k-1}),
\end{align*}
with $\mu = \sup_{t \in \RR_+} \EE(X_t) < \infty,$ and that similarly, \[
\EE(Y_{k+1}^2) \leq \sigma^2 \mu (t_{i^* + k} - t_{i^*+ k-1}) = \sigma^2 \mu 2^{(i^*+ k-2)}.
\] 
All in all, with Lemma \ref{kole}, we can estimate the first term in \eqref{eq:disc_cont_decomp} by
\begin{align*}
\frac{4}{\vare^2}& \left(3 \sigma^2 \mu \sum_{k=1}^{n-i^*+1} 4^{-(i^* +k -1)} 2^{(i^* +k-1)} + 2 \sigma^2 \mu \sum_{k=1}^{n-i*+1} 4^{-(i^* +k - 1)} (2^{(i^* +k -2) + (i^* +k-1)})^{1/2} \right. \\
  & \left. \quad +2 \sigma^2 \mu \sum_{k=1}^{n-i*+1} 4^{-(i^* + k - 1)} 2^{(i^* + k -2)} \right) \\
  &\leq \frac{4}{\vare^2} \left( 3 \sigma^2 \mu 2^{-i^*} \sum_{k=1}^{\infty}2^{-(k-1)} + 2 \sigma^2 \mu 2^{-i^*} \sum_{k=1}^{\infty}2^{-(k-1)} + 2\sigma^2 \mu 2^{-i^*} \sum_{k=1}^{\infty}2^{-(k-1)}\right) = \frac{56}{\vare^2}\sigma^2 \mu 2^{-i^*}.
\end{align*}
This does not depend on $n$ (hence, on $T$), and since $i^* \to \infty$ as $K \to \infty$, we have that the first term in \eqref{eq:disc_cont_decomp} converges to zero as $\rho T \to \infty$ and then $K \to \infty$.

For the second term in \eqref{eq:disc_cont_decomp} we will use Doob's submartingale inequality \citep[see, e.g., ][Theorem 1.3.8. (i)]{Karatzas_91} to the submartingales 
\[
N_{t,i} := (M_{t_{i+1}+t} - M_{t_{i+1}})^2, \quad t \in [0, t_i - t_{i+1}], \quad i = i^*,\ldots,n,
\]
for which clearly
\[
\PP\left(\sup_{t_{i+1} < t < t_{i}} \left| (\rho T - t_{i})^{-1} (M_t - M_{t_{i+1}}) \right| > \frac{\vare}{2} \right) = \PP\left(\sup_{0 \leq t \leq t_i - t_{i+1}} N_{t,i} > \frac{\vare^2(\rho T-t_i)^2}{4}\right).
\]
The inequality states that
\[
\PP\left(\sup_{0 \leq t \leq t_i - t_{i+1}} N_{t,i} > \frac{\vare^2(\rho T-t_i)^2}{4}\right) \leq \frac{4 \EE(N_{t_i - t_{i+1}})}{\vare^2 (\rho T-t_i)^2} = \frac{4 \EE(M_{t_i} - M_{t_{i+1}})^2}{\vare^2 (\rho T-t_i)^2} \leq \frac{4\sigma^2 \mu (t_i - t_{i+1})}{\vare^2 (\rho T-t_i)^2}.
\]
Now, in our present setting, \[
t_i - t_{i+1} \leq (\rho T - 2^{i-1}) - (\rho T - 2^{i}) = 2^{i-1} \quad \text{and} \quad (\rho T - t_i)^2 \geq 2^{2i-4}.
\]
Thus, the second term in \eqref{eq:disc_cont_decomp} can be estimated from above by
\[
\frac{\sigma^2 \mu \vare^2}{4}\sum_{i=i^*}^{n} 2^{-i+3} \leq \frac{\sigma^2 \mu \vare^2}{4} 2^{-i^*+3} \sum_{i=0}^{\infty} 2^{-i}.
\]
Again, clearly this does not depend on $n$ (thus, $T$) and converges to zero as $i^* \to \infty$ (and thus, as $K \to \infty$). 
This suffices for our statement. \proofend

\section*{Acknowledgement}
The authors are grateful to Professor P\'eter Major at the University of Szeged for supplying the basic idea of the proof of Lemma \ref{lem:M_sup}.

\bibliographystyle{apalike}
\bibliography{master}

\indent
      {\sc Gyula Pap},
      Bolyai Institute, University of Szeged,
      Aradi v\'ertan\'uk tere 1, H--6720 Szeged, Hungary.
      E--mail: papgy@math.u-szeged.hu

\indent
      {\sc Tam\'as T. Szab\'o},
      Bolyai Institute, University of Szeged,
      Aradi v\'ertan\'uk tere 1, H--6720 Szeged, Hungary.
      Tel.: +36-62-343882, Fax: +36-62-544548,
      E--mail: tszabo@math.u-szeged.hu

\end{document}